\documentclass[11pt]{amsart}
\usepackage{}

\usepackage{amsmath}
\allowdisplaybreaks
\usepackage{amsfonts}
\usepackage{amssymb}
\usepackage[all]{xy}           

\usepackage{bbding}
\usepackage{txfonts}
\usepackage{amscd}

\usepackage[shortlabels]{enumitem}
\usepackage{ifpdf}
\ifpdf
  \usepackage[colorlinks,final,backref=page,hyperindex]{hyperref}
\else
  \usepackage[colorlinks,final,backref=page,hyperindex,hypertex]{hyperref}
\fi
\usepackage{tikz}
\usepackage[active]{srcltx}

\makeatletter

\newtheorem{df}{Definition}[section]
\newtheorem{thm}{Theorem}[section]

\newtheorem{rem}{Remark}[section]

\newtheorem{prop}{Proposition}[section]
\newtheorem{exa}{Example}[section]
\newtheorem{lem}{Lemma}[section]

\numberwithin{equation}{section}

\begin{document}

\date{}
\title{On HOM-NAMBU POISSON SUPERALGEBRAS STRUCTURES AND THEIR REPRESENTATIONS}
\author{OTHMEN NCIB }
\address{University of Gafsa, Faculty of Sciences Gafsa, 2112 Gafsa, Tunisia}
\email{othmenncib@yahoo.fr}

 \maketitle{}

 \begin{abstract} The purpose of this paper is to study some results of constructions on Hom-Poisson superalgebras we use the representations and Rota-Baxter operators. We introduce the structures of $n$-ary Hom-Nambu Poisson superalgebras and their representations and we study relationships between Hom-Poisson superalgebras and its induced $n$-ary Hom-Nambu Poisson superalgebras. The same concept is given for the representations of $n$-ary Hom-Nambu Poisson superalgebras.

 \end{abstract}
 
 \textbf{key words:} Hom-Lie superalgebras, $n$-ary Hom-Nambu superalgebras, Hom-Poisson superalgebras, $n$-ary Hom-Nambu-Poisson superalgebras, Representation, Rota-Baxter operators. \\
 
 
\tableofcontents

\section*{Introduction}

The notion of Hom-algebras is first appeared in $q$-deformations of algebras of vector fields in physics contexts, These structures include their classical counterparts and open more possibilities for deformations, Hom-algebra extensions of cohomological structures and representations, formal deformations of Hom-associative and Hom-Lie algebras, Hom-Lie admissible Hom-coalgebras, Hom-coalgebras, Hom-Hopf algebras \cite{Ammar-Ejbehi-Makhlouf,Benayadi-Makhlouf,Larsson-Silvestrov,MakhloufSilvestrov1,MakhloufSilvestrov2,MakhloufSilvestrov3,Sheng,Yau1,Yau2}. The most useful algebraic structures of this type are the Home associative algebras which is twisted version of the associative algebras introduced by Makhlouf and Silvestrov in \cite{Makhlouf-Silvestrov} and the Hom-Lie algebras and more general quasi-Hom-Lie algebras structures were introduced first by Hartwig, Larsson and Silvestrov in \cite{Hartwig-Larsson-Silvestrov}, where the general quasi-deformations and discretizations of Lie algebras of vector fields using more general $\sigma$-derivations and a general method for construction of deformations of Witt and Virasoro type algebras based on twisted derivations have been developed.\\

In physics, algebras with $n$-ary compositions appear in many different contexts. Let us mention a few of them. Ternary
algebras can be used to construct solutions of the Yang–Baxter equation \cite{Okubo}, which first appeared in statistical
mechanics \cite{R-J-Baxter1,R-J-Baxter2}. Nambu mechanics \cite{Y-Nambu} involves an n-ary product that satisfies the n-ary Nambu identity, which is an
n-ary generalization of the Jacobi identity . Bagger–Lambert algebras \cite{Bagger-Lambert} are ternary Nambu algebras
with some extra structures, and they appear in the study of string theory and M-branes. Ternary algebras are used in \cite{Gunaydin,Gunaydin-Hyun}
to construct superconformal algebras.
Generalizations of $n$-ary Nambu and Nambu–Lie algebras, called $n$-ary Hom–Nambu and Hom–Nambu–Lie algebras, were
introduced in \cite{Ataguema-Makhlouf-Silvestrov} by Ataguema, Makhlouf, and Silvestrov. In these Hom-type algebras, the $n$-ary Nambu identity is relaxed
using $n-1$ linear maps, called the twisting maps, resulting in the $n$-ary Hom–Nambu identity . When these
twisting maps are all equal to the identity map, one recovers $n$-ary Nambu and Nambu–Lie algebras.\\

Nambu-Poisson algebras structures were intrioduced in \cite{Takhtajan-L} such as a generalization of Poisson-algebras, arose
from the study of generalized Nambu mechanics and its quantization, which plays an important role in physics, mainly in string theory. Nambu–Poisson algebras in noncommutative
geometry and the quantum geometry of branes in M-theory have been considered in \cite{DeBellis-Samann-Szabo1} and the theory of Nambu–Poisson algebras provides an useful method to describe algebraical and geometrical structures of Nambu–Poisson manifolds \cite{Das,DeBellis-Samann-Szabo2,Song-Jiang,Vallejo}. The Hom type of this notion (i.e., $\alpha$-twisted) is the structure of Nambu-poisson algebra endowed with a homomorphism defined algebraically as an $n$-Hom-Lie algebra  $(\mathcal{N},\{,\},\alpha)$ endowed with a Hom-associative product $\mu$ and with both products compatibel via a Hom–Leibnize identity
\begin{equation}\label{n-Leibniz-identity-intr}
\{\alpha(x_1),\cdots,\alpha(x_{n-1}),\mu(y,z)\}=\mu(\{x_1,\cdots,x_{n-1},y\},\alpha(z))+\mu(\{x_1,\cdots,x_{n-1},z\},\alpha(y)),
\end{equation}
which is gives by many authors (see \cite{AmriMakhlouf,Laraiedh-Silvestrov,MakhloufSilvestrov3,SHANSHAN-MAKHLOUF-LINASONG,Yau3}). More applications of (Hom)-Lie superalgebras, (Hom)-Lie color algebras and some results on the structures of (Hom)-Poisson superalgebras have been widely investigated \cite{Ammar-Ayadi-Mabrouk-Makhlouf,Ammar-Makhlouf,Ammar-Makhlouf-Saadoui,Liu-Hu,Ben Hassine-Mabrouk-Ncib,THOMAS,IVAN,Chunyue-Qingcheng-Zhu}. Some classifications of Hom-associative superalgebras, Hom-Lie superalgebras and Hom-Poisson superalgebras are studied in \cite{Qingcheng-Chunyue-Wei1}.\\

Representations and Rota-Baxter operators of different algebraic structures (seeN\cite{})
are an important subject of study in algebra and diverse areas. They appear
in many fields of mathematics and physics. In particular, they appear in deformation and
 cohomology  theory among other areas. They have
various applications relating algebra to geometry and allow the construction of new algebraic structures. The aim of this paper is to study $n$-ary Hom-Nambu Poisson superalgebras and representation over $n$-ary Hom-Nambu Poisson superalgebras.\\

The paper is organized as follows. In section $1$, we give the definition and some important constructions of Hom-associative superalgebras, Hom-Lie superalgebras and Hom-Poisson superalgebnras. We
also recall some examples for these structures and we give some results coserning the representations and Rota-Baxter concerning these structures. In section $2$, we define the notion of $n$-ary Hom-Nambu Poisson superalgebras and we provide a construction for a multiplicative $n$-ary Hom-Nambu Poisson superalgebras using multiplicative Hom-Poisson superalgebras. In section $3$ the notion of representation of $n$-ary
Hom-Nambu Poisson superalgebras are introduced and some results obtained.\\

Throughout this paper, $\mathbb{K}$ denotes a field of characteristic zero. All vector spaces and
algebras are $\mathbb{Z}_2$-graded over $\mathbb{K}$. Such a superspace over $\mathbb{K}$ is a $\mathbb{Z}_2$ graded $\mathbb{K}$-linear space $V =V_0\oplus V_1$. The
element of $V_k,\; k \in \{0, 1\}$, are said to be homogenous and of parity $k$. We denoted by $|x|$, the parity of a homogeneous element $x$.

 \section{Preliminaries and some results }

In this section we introduce the notions of Hom-associative superalgebras, Hom-Lie superalgebras and Hom-Poisson superalgebras in which we give some results and some examples.
\begin{df}
A Hom-associative superalgebre is a triplet $(\mathcal{A},\mu,\alpha)$ consisting of a $\mathbb{Z}_2$-graded vector space $\mathcal{A}$, an even linear map
 $\alpha:\mathcal{A}\rightarrow\mathcal{A}$, (i.e: $\alpha(\mathcal{A}_i)\subseteq\mathcal{A}_i)$ and an even bilinear map
 $\mu:\mathcal{A}\rightarrow\mathcal{A}$, (i.e: $\mu(\mathcal{A}_i,\mathcal{A}_j)\subseteq\mathcal{A}_{i+j})$ satisfying the following condition
 for all $x,y,z\in\mathcal{H}(\mathcal{A})$:
 $$ass_{\mu}(x,y,z)=0\;(\text{Hom-associativity}),$$
 where $ass_{\mu}(x,y,z)=\mu(\mu(x,y),\alpha(z))-\mu(\alpha(x),\mu(y,z)),\;\forall x,y,z\in\mathcal{H}(\mathcal{A}).$
\end{df}
If in addition $\mu$ is super-commutative (i.e: $\mu(x,y)=(-1)^{|x||y|}\mu(y,x)\;\forall\;x,y\in\mathcal{H}(\mathcal{A})$),
the Hom-associative superalgebra $(\mathcal{A},\mu,\alpha)$ is said to be commutative Hom-associative superalgebra. .
\begin{exa}
Let $\mathcal{A}=\mathcal{A}_{\overline{0}}\oplus \mathcal{A}_{\overline{1}}$ be a $3$-dimensional $\mathbb{Z}_2$-
graded vector space, where $\mathcal{A}_{\overline{0}}=<e_1,e_2>$
 and $\mathcal{A}_{\overline{1}}=<e_3>$. The nonzero products given by
 $$\mu(e_1,e_2)=-e_1,\;\mu(e_2,e_2)=e_1+e_2,$$
 and the even linear map $\alpha:\mathcal{A}\rightarrow\mathcal{A}$ defined on the basis of $\mathcal{A}$ by
 $$\alpha(e_1)=-e_1,\;\alpha(e_2)=e_1+e_2\;\text{and}\;\alpha(e_3)=0,$$ defined on $\mathcal{A}$ a structure of Hom-associative superalgebra.
\end{exa}
\begin{df}
A Hom-Lie superalgebra is a triple $(\mathcal{A}, [~~,~~],\alpha)$ consisting of a $\mathbb{Z}_2$-graded vector space $\mathcal{A}$, an even
bilinear map $[~~,~~] : \mathcal{A}\otimes \mathcal{A} \longrightarrow \mathcal{A},~~([\mathcal{A}_i,\mathcal{A}_j]\subseteq \mathcal{A}_{i+j},
~~\forall~~i,j\in \mathbb{Z}_2)$ and an even linear map  $\alpha:\mathcal{A}\rightarrow\mathcal{A}$  satisfying:
\begin{eqnarray}
 \label{H-skwesym}
 &&[x,y] = -(-1)^{|x||y|}[y,x]\quad \text{( super-skew-symmetry)},\\
\label{H-sJ}
&&  [\alpha(x),[y,z]]=[[x,y],\alpha(z)]+(-1)^{|x||y|}[\alpha(y),[x,z]] \quad \text{(Hom super-Jacobi identity)},
\end{eqnarray}
$\forall\ x,y,z \in \mathcal{H}(\mathcal{A})$.
\end{df}
\begin{exa}
Let $\mathcal{A}=\mathcal{A}_{\overline{0}}\oplus \mathcal{A}_{\overline{1}}$ be a $2$-dimensional $\mathbb{Z}_2$-
graded vector space, where $\mathcal{A}_{\overline{0}}=<e_1>$
 and $\mathcal{A}_{\overline{1}}=<e_2>$. We consider the nonzero product given by
$$[e_2,e_2]=2\lambda e_1,$$
and the even linear map $\alpha:\mathcal{A}\rightarrow\mathcal{A}$ defined by $$\alpha(e_1)=\lambda e_1\;\text{and}\;\alpha(e_2)=\lambda e_2,$$
where $\lambda\in\mathbb{C}$. The triplet $(\mathcal{A},[\cdot,\cdot],\alpha)$ is a Hom-Lie superalgebra.
\end{exa}
\begin{prop}(\cite{Ammar-Makhlouf})\label{Hom-Lie-super-twist}
Let $(\mathcal{A},\mu,\alpha)$ be a Hmo-associative superalgebra, then $(\mathcal{A},[\cdot,\cdot]_\mu,\alpha)$ is a Hom-Lie superalgebra, where
$$[x,y]_\mu=\mu(x,y)-(-1)^{|x||y|}\mu(y,x),$$
for all $x,y\in\mathcal{H}(\mathcal{A})$.
\end{prop}
\begin{df}\textbf{(Non commutative Hom-Poisson superalgebras):}
A non-commutative Hom-Poisson superalgebra is a quadruple $(\mathcal{A}, [\cdot, \cdot], \mu, \alpha)$ consists of
\begin{enumerate}
\item A Hom-Lie superalgebra $(\mathcal{A}, [\cdot, \cdot],\alpha)$,
\item A Hom-associative superalgebra $(\mathcal{A},\mu, \alpha)$,
\end{enumerate}
such that the Hom-Leibniz superalgebra identity
\begin{equation}\label{Homleibnizident}
[\alpha(x),\mu(y,z)]=\mu([x,y],\alpha(z))+(-1)^{|x||y|}\mu(\alpha(y),[x,z])
\end{equation}
is satisfied for all $x,y,z\in\mathcal{H}(\mathcal{A})$. A Hom-Poisson superalgebra is a non-commutative Hom-Poisson superalgebra
 $(\mathcal{A}, [\cdot, \cdot], \mu, \alpha)$ in
which $\mu$ is super-commutative.\\
\begin{rem}
If $(\mathcal{A}, [\cdot, \cdot], \mu, \alpha)$ is a Hom-Poisson superalgebra then the identity \eqref{Homleibnizident} is equivalent to
\begin{equation}\label{equiv-Homleibnizident}
[\alpha(x),\mu(y,z)]=\mu([x,y],\alpha(z))+(-1)^{|y||z|}\mu([x,z],\alpha(y)).
\end{equation}
\end{rem}
If $\alpha$ is furthermore a Poisson automorphism, that is, a linear bijective on such that
$\alpha([x, y]) = [\alpha(x), \alpha(y)]$ and $\alpha(\mu(x,y)) = \mu(\alpha(x),\alpha(y))$ for any $x, y\in\mathcal{H}(\mathcal{A})$,
 then $\mathcal{A}$ is called a regular Hom-Poisson superalgebra.

\end{df}

\begin{exa}
Let $\mathcal{A} = \mathcal{A}_{\overline{0}} \oplus \mathcal{A}_{\overline{1}}$ be a $3$-dimensional superspace, where
$\mathcal{A}_{\overline{0}}=<e_1,e_2>$ and $\mathcal{A}_{\overline{1}}=<e_3>$. We define on the basis of $\mathcal{A}$, the following linear map
$$\alpha(e_1)=e_1\;,\;\alpha(e_2)=e_1+e_2,$$
and the following nonzero products
$$\mu(e_1,e_2)=e_1\;,\;\mu(e_2,e_2)=e_1+e_2,$$
$$[e_1,e_2]=e_1.$$
Then the quadruplet $(\mathcal{A},[\cdot,\cdot],\mu,\alpha)$ is a non-commutative Hom-Poisson superalgebra.
\end{exa}
 \begin{prop}
Let $(\mathcal{A}, \mu, \alpha)$ be a Hom-associative superalgebra. Then
$$\mathcal{A}_\mu=(\mathcal{A}, [\cdot, \cdot]_\mu, \mu, \alpha)$$
is a non-commutative Hom-Poisson superalgebra, where
$$[x,y ]_\mu=\mu(x,y)-(-1)^{|x||y|}\mu(y,x),\;\forall x,y\in \mathcal{H}(\mathcal{A}).$$

\end{prop}

\begin{proof}
By Proposition \ref{Hom-Lie-super-twist}, $(\mathcal{A},[\cdot,\cdot]_\mu,\alpha)$ is a Hom-Lie superalgebra.\\

Let $x,y,z\in\mathcal{H}(\mathcal{A})$, we have
\begin{align*}
\mu([x,y]_\mu,\alpha(z))+(-1)^{|x||y|}\mu(\alpha(y),[x,z]_\mu)&=\mu(\mu(x,y),\alpha(z))-(-1)^{|x||y|}\mu(\mu(y,x),\alpha(z))+(-1)^{|x||y|}\mu(\alpha(y),\mu(x,z))\\&-(-1)^{|x|(|y|+|z|)}\mu(\alpha(y),\mu(z,x))\\&=\mu(\alpha(x),\mu(y,z))-(-1)^{|x||y|}\mu(\mu(y,x),\alpha(z))+(-1)^{|x||y|}\mu(\mu(y,x),\alpha(z))\\&-(-1)^{|x|(|y|+|z|)}\mu(\mu(y,z),\alpha(x)))\\&=\mu(\alpha(x),\mu(y,z))-(-1)^{|x|(|y|+|z|)}\mu(\mu(y,z),\alpha(x)))\\&=[\alpha(x),\mu(y,z)]_\mu.
\end{align*}
This shows that the Hom-Leibniz superalgebra identity is satisfied, which gives that $(\mathcal{A},[\cdot,\cdot]_\mu,\mu,\alpha)$ is a Hom-Poisson superalgebra. 
\end{proof}

\begin{df}
Let $(\mathcal{A}, [\cdot, \cdot], \mu, \alpha)$ be Hom-Poisson superalgebras.\\
\begin{enumerate}
\item $\mathcal{A}$ is multiplicative if
$$\alpha\circ[\cdot, \cdot]=[\cdot, \cdot]\circ\alpha^{\otimes2}~~~~and~~~~\alpha\circ\mu=\mu\circ\alpha^{\otimes2}.$$
\item Let $(\mathcal{B}, [\cdot, \cdot]', \mu', \beta)$ be another Hom-Poisson superalgebra. A weak morphism
$f : \mathcal{A}\rightarrow \mathcal{B}$ is an even linear
map such that $$f\circ[\cdot, \cdot]=[\cdot, \cdot]'\circ f^{\otimes2}~~~~and~~~~f\circ\mu=\mu'\circ f^{\otimes2}.$$
\item A morphism $f : \mathcal{A} \rightarrow \mathcal{B}$ is a weak morphism such that $f\circ\alpha=\beta\circ f$.
\end{enumerate}

\end{df}

 \subsection{Derivation and Rota-Baxter operator of Hom-Poisson Superalgebras. }

\begin{df}
\begin{enumerate}
\item A derivation of a Hom-associative superalgebra $(\mathcal{A},\mu,\alpha)$ is a linear map $\mathfrak{D}:\mathcal{A}\rightarrow\mathcal{A}$ satisfying
$$\mathfrak{D}\alpha=\alpha\mathfrak{D},$$
$$\mathfrak{D}(\mu(x,y))=\mu(\mathfrak{D}(x),y)+(-1)^{|\mathfrak{D}||x|}\mu(x,\mathfrak{D}(y)),\;\forall x,y\in\mathcal{H}(\mathcal{A}).$$
\item A derivation of a Hom-Lie superalgebra $(\mathcal{A},[\cdot,\cdot],\alpha)$ is a linear map $\mathfrak{D}:\mathcal{A}\rightarrow\mathcal{A}$ satisfying
$$\mathfrak{D}\alpha=\alpha\mathfrak{D},$$
$$\mathfrak{D}([x,y])=[\mathfrak{D}(x),y]+(-1)^{|\mathfrak{D}||x|}[x,\mathfrak{D}(y)],\;\forall x,y\in\mathcal{H}(\mathcal{A}).$$
\item A derivation of a Hom-Poisson superalgebra $(\mathcal{A},[\cdot,\cdot],\mu,\alpha)$ is a linear map $\mathfrak{D}:\mathcal{A}\rightarrow\mathcal{A}$ satisfying
\begin{eqnarray}
\mathfrak{D}\alpha&=&\alpha\mathfrak{D}\label{Poisson-Deriv1},\\
\mathfrak{D}(\mu(x,y))&=&\mu(\mathfrak{D}(x),y)+(-1)^{|\mathfrak{D}||x|}\mu(x,\mathfrak{D}(y))\label{Poisson-Deriv2},\\
\mathfrak{D}([x,y])&=&[\mathfrak{D}(x),y]+(-1)^{|\mathfrak{D}||x|}[x,\mathfrak{D}(y)]\label{Poisson-Deriv3},\;\forall x,y\in\mathcal{H}(\mathcal{A}).
\end{eqnarray}
\end{enumerate}
\end{df}
 For any $x\in\mathcal{A}$, define $ad_x\in End_\mathbb{K}(\mathcal{A})$ by $ad_x(y)=[x,y]$, for any $y\in\mathcal{A}$.
 Then the Hom super-Jacobi identity \eqref{H-sJ} can be written as
 \begin{equation}\label{ident-Jac-adj}
  ad_{[x,y]}(\alpha(z))=ad_{\alpha(x)}\circ ad_y(z)-(-1)^{|x||y|}ad_{\alpha(y)}\circ ad_x(z)
 \end{equation}
 for all $x, y, z\in\mathcal{H}(\mathcal{A})$.
 \begin{rem}
 The Hom-Leibniz identity \eqref{Homleibnizident} says that $ad_x$ is a derivation with respect to the Hom-associative product, i.e. the equation
  \eqref{Homleibnizident} rewrite in the form of the following expression
 $$ad_{\alpha(x)}(\mu(y,z))=\mu(ad_x(y),\alpha(z))+(-1)^{|x||y|}\mu(\alpha(y),ad_x(z)),$$
 for all $x,y,z\in\mathcal{H}(\mathcal{A})$.
 \end{rem}
A Rota-Baxter operator of weight $\lambda\in\mathbb{K}$ on a Hom-associative superalgebra $(\mathcal{A},\mu,\alpha)$ is an even linear map
 $\mathcal{R}:\mathcal{A}\to\mathcal{A}$ commuting with $\alpha$ and satisfying
\begin{equation}\label{rota-baxter-Hom-assoc}
    \mu(\mathcal{R}(x),\mathcal{R}(y))=\mathcal{R}\Big(\mu(\mathcal{R}(x),y)+\mu(x,\mathcal{R}(y))+\lambda\mu(x,y)\Big),\forall x,y\in\mathcal{H}(\mathcal{A}).
\end{equation}

In \cite{Abdaoui-Mabrouk-Makhlouf}, the authors are defined a Rota-Baxter operator of weight $\lambda\in\mathbb{K}$ on a Hom-Lie superalgebras
$(\mathcal{A},[\cdot,\cdot],\alpha)$ as an even linear map $\mathcal{R}:\mathcal{A}\to\mathcal{A}$ commuting with
$\alpha$ and satisfies the following condition
\begin{equation}\label{hom-superalg-Rota-Baxter}
  [\mathcal{R}(x),\mathcal{R}(y)]=\mathcal{R}\Big([\mathcal{R},y]+[x,\mathcal{R}(y)]+\lambda[x,y]\Big),
\end{equation}
for all $x,y\in\mathcal{H}(\mathcal{A}).$
\begin{prop}
Let $R$ be a Rota-Baxter operator of weight zero on a Hom-Lie superalgebra $(\mathcal{A},[\cdot,\cdot],\alpha)$ (resp. Hom-associative superalgebra $(\mathcal{A},\mu,\alpha)$). Then, $R$ is bijective if and only if $R^{-1}$ is a derivation of $(\mathcal{A},[\cdot,\cdot],\alpha)$ (resp. derivation of $(\mathcal{A},\mu,\alpha)$).
\end{prop}
\begin{df}
A Rota-Baxter operator of Weight $(\lambda,\lambda')$ on a Hom-Poisson superalgebra $(\mathcal{A},[\cdot,\cdot],\mu,\alpha)$
 is a pair $(\mathcal{R},\mathcal{R}')$ consisting of an even commuting linear maps on $\mathcal{A}$ in which $\mathcal{R}$
 is a Rota-Baxter operator of weight $\lambda$ on the Hom-associative superalgebra $(\mathcal{A},\mu,\alpha)$ and
 $\mathcal{R}'$ is a Rota-Baxter operator of weight $\lambda'$ on the Hom-Lie superalgebra $(\mathcal{A},[\cdot,\cdot],\alpha)$.
\end{df}
\begin{lem}\label{Rota-Baxt_assoc-superalg}
Let $\mathcal{R}$ be a Rota-Baxter operator of weight $\lambda$ on a Hom-associative superalgebra $(\mathcal{A},\mu,\alpha)$. Then the bilniear product
\begin{equation}\label{struc-assoc-Rot-Bax}
  \mu_{\mathcal{R}}(x,y)= \mu(\mathcal{R}(x),y)+\mu(x,\mathcal{R}(y))+\lambda\mu(x,y),\forall x,y\in\mathcal{H}(\mathcal{A}),
\end{equation}
define a structure of a Hom-associative superalgebra on $\mathcal{A}$.
\end{lem}

\begin{proof}
By eq. \eqref{rota-baxter-Hom-assoc}, we clearly notice that $\mathcal{R}(\mu_{\mathcal{R}}(x,y))=
\mu(\mathcal{R}(x),\mathcal{R}(y)),\forall x,y\in\mathcal{H}(\mathcal{A})$.\\
Let $x,y,z\in\mathcal{H}(\mathcal{A})$, we have
\begin{align*}
\mu_{\mathcal{R}}\Big(\mu_{\mathcal{R}}(x,y),\alpha(z)\Big)&=\mu\Big(\mathcal{R}(\mu_{\mathcal{R}}(x,y)),\alpha(z)\Big)+
\mu\Big(\mu_{\mathcal{R}}(x,y),\mathcal{R}(\alpha(z))\Big)+\lambda\mu\big(\mu_{\mathcal{R}}(x,y),\alpha(z)\Big)\\&=
\mu\Big(\mu(\mathcal{R}(x),\mathcal{R}(y)),\alpha(z)\Big)+\mu\Big(\mu(\mathcal{R}(x),y),\alpha(\mathcal{R}(z))\Big)\\&
+\mu\Big(\mu(x,\mathcal{R}(y)),\alpha(\mathcal{R}(z))\Big)+\lambda\mu\Big(\mu(x,y),\alpha(\mathcal{R}(z))\Big)\\&
+\lambda\mu\Big(\mu(\mathcal{R}(x),y),\alpha(z)\Big)+\lambda\mu\Big(\mu(x,(y)),\alpha(z)\Big)+\lambda^2\mu\Big(\mu(x,y),\alpha(z)\Big).
\end{align*}
Similarly, we have
\begin{align*}
\mu_{\mathcal{R}}\Big(\alpha(x),\mu_{\mathcal{R}}(y,z)\Big)&=\mu\Big(\alpha(\mathcal{R}(x)),\mu(\mathcal{R}(y),z)\Big)+
\mu\Big(\alpha(\mathcal{R}(x)),\mu(y,\mathcal{R}(z))\Big)\\&+\lambda\mu\Big(\alpha(\mathcal{R}(x)),\mu(y,z)\Big)+
\mu\Big(\alpha(x),\mu(\mathcal{R}(y),\mathcal{R}(z)))\Big)\\&+\lambda\mu\Big(\alpha(x),\mu(\mathcal{R}(y),z))\Big)+
\lambda\mu\Big(\alpha(x),\mu(y,\mathcal{R}(z)))\Big)+\lambda^2\mu\Big(\alpha(x),\mu(y,z))\Big).
\end{align*}
According to eq. \eqref{rota-baxter-Hom-assoc} and by a direct computation, we get
$$\mu_{\mathcal{R}}\Big(\mu_{\mathcal{R}}(x,y),\alpha(z)\Big)=\mu_{\mathcal{R}}\Big(\alpha(x),\mu_{\mathcal{R}}(y,z)\Big),\forall x,y,z\in\mathcal{H}(\mathcal{A}).$$
Then $(\mathcal{A},\mu_{\mathcal{R}},\alpha)$ is a Hom-associative superalgebra.
\end{proof}

\begin{lem}\label{Lie-Rota-Baxt_superalg}
Let $(\mathcal{A},[\cdot,\cdot],\alpha)$ be a Hom-Lie superalgebra and $\mathcal{R}:\mathcal{A}\to\mathcal{A}$
 a Rota-Baxter operator of weight $\lambda$ on $\mathcal{A}$. Then $(\mathcal{A},[\cdot,\cdot]_{\mathcal{R}},\alpha)$ is a Hom-Lie superalgebra, where
\begin{equation}\label{struc-Lie-Rot-Bax}
 [x,y]_{\mathcal{R}}=[\mathcal{R}(x),y]+[x,\mathcal{R}(y)]+\lambda[x,y],\forall x,y \in\mathcal{H}(\mathcal{A}).
\end{equation}
\end{lem}
\begin{proof}
It's easy to see that $[\cdot,\cdot]_{\mathcal{R}}$ is super-skew-symmetric.\\
Let $x,y,z\in\mathcal{H}(\mathcal{A})$, we have
\begin{align*}
    [\alpha(x),[y,z]_{\mathcal{R}}]_{\mathcal{R}}&=[\alpha(\mathcal{R}(x)),[y,z]_{\mathcal{R}}]+
    [\alpha(x),\mathcal{R}[y,z]_{\mathcal{R}}]+\lambda[\alpha(x),[y,z]_{\mathcal{R}}]\\&=
    [\alpha(\mathcal{R}(x)),[\mathcal{R}(y),z]]+[\alpha(R(x)),[y,\mathcal{R}(z)]]+
    \lambda[\alpha(\mathcal{R}(x)),[y,z]]\\&+[\alpha(x),[\mathcal{R}(y),\mathcal{R}(z)]]+
    \lambda[\alpha(x),[\mathcal{R}(y),z]]+\lambda[\alpha(x),[y,\mathcal{R}(z)]]+\lambda^2[\alpha(x),[y,z]].
\end{align*}
In the same way, we have
\begin{align*}
    [\alpha(y),[z,x]_{\mathcal{R}}]_{\mathcal{R}}&=[\alpha(\mathcal{R}(y)),[\mathcal{R}(z),x]]+
    [\alpha(R(y)),[z,\mathcal{R}(x)]]+\lambda[\alpha(\mathcal{R}(y)),[z,x]]\\&+
    [\alpha(y),[\mathcal{R}(z),\mathcal{R}(x)]]+\lambda[\alpha(y),[\mathcal{R}(z),x]]+\lambda[\alpha(y),[z,\mathcal{R}(x)]]+
    \lambda^2[\alpha(y),[z,x]],
\end{align*}
and
\begin{align*}
    [\alpha(x),[y,z]_{\mathcal{R}}]_{\mathcal{R}}&=[\alpha(\mathcal{R}(z)),[\mathcal{R}(x),y]]+
    [\alpha(R(z)),[x,\mathcal{R}(y)]]+\lambda[\alpha(\mathcal{R}(z)),[x,y]]\\&+[\alpha(z),[\mathcal{R}(x),\mathcal{R}(y)]]+
    \lambda[\alpha(z),[\mathcal{R}(x),y]]+\lambda[\alpha(z),[x,\mathcal{R}(y)]]+\lambda^2[\alpha(z),[x,y]].
\end{align*}
We can easly conclude that the identity \eqref{H-sJ} is holds on $[\cdot,\cdot]_{\mathcal{R}}$, which gives that
$(\mathcal{A},[\cdot,\cdot]_{\mathcal{R}},\alpha)$ is a Hom-Lie superalgebras.
\end{proof}
\begin{thm}
Let $(\mathcal{R},\mathcal{R})$ is a Rota-Baxter operator of weight $(\lambda,\lambda)$ on a Hom-Poisson superalgebra
$(\mathcal{A},[\cdot,\cdot],\mu,\alpha)$. Then $(\mathcal{A},[\cdot,\cdot]_{\mathcal{R}},\mu_{\mathcal{R}},\alpha)$ is a Hom-Poisson superalgebra.
\end{thm}
\begin{proof}
By Lemmas \ref{Rota-Baxt_assoc-superalg} and \ref{Lie-Rota-Baxt_superalg}, we deduce that
$(\mathcal{A},\mu_{\mathcal{R}},\alpha)$ is a Hom-associative superalgebra and
$(\mathcal{A},[\cdot,\cdot]_{\mathcal{R}},\alpha)$ is a Hom-Lie superalgebra.\\

Let $x,y,z\in\mathcal{H}(\mathcal{A})$.\\
\begin{align*}
    \bullet)\;[\alpha(x),\mu_{\mathcal{R}}(y,z)]_{\mathcal{R}}&=[\alpha({\mathcal{R}}(x)),\mu_{\mathcal{R}}(y,z)]+
    [\alpha(x),{\mathcal{R}}(\mu_{\mathcal{R}}(y,z))]+\lambda[\alpha(x),\mu_{\mathcal{R}}(y,z)]\\&=
    [\alpha({\mathcal{R}}(x)),\mu(\mathcal{R}(y),z)]+[\alpha({\mathcal{R}}(x)),\mu(y,\mathcal{R}(z))]+
    \lambda[\alpha({\mathcal{R}}(x)),\mu(y,z)]\\&+[\alpha(x),\mu(\mathcal{R}(y),\mathcal{R}(z))]+
    \lambda[\alpha(x),\mu(\mathcal{R}(y),z)]+\lambda[\alpha(x),\mu(y,\mathcal{R}(z))]+\lambda^2[\alpha(x),\mu(y,z)].
\end{align*}
\begin{align*}
\bullet)\;\mu_{\mathcal{R}}([x,y]_{\mathcal{R}},\alpha(z)) &=\mu(\mathcal{R}[x,y]_{\mathcal{R}},\alpha(z))+
\mu([x,y]_{\mathcal{R}},\alpha(\mathcal{R}(z)))+\lambda\mu([x,y]_{\mathcal{R}},\alpha(z))\\&=
\mu([\mathcal{R}(x),\mathcal{R}(y)],\alpha(z))+ \mu([\mathcal{R}(x),y],\alpha(\mathcal{R}(z)))+
\mu([x,\mathcal{R}(y)],\alpha(\mathcal{R}(z)))\\&+\lambda\mu([x,y],\alpha(\mathcal{R}(z)))+
\lambda\mu([\mathcal{R}(x),y],\alpha(z))+\lambda\mu([x,\mathcal{R}(y)],\alpha(z))+\lambda^2\mu([x,y],\alpha(z)).
\end{align*}
Similarly, we have
\begin{align*}
\bullet)\;\mu_{\mathcal{R}}(\alpha(y),[x,z]_{\mathcal{R}})&=\mu(\alpha(\mathcal{R}(y)),[\mathcal{R}(x),z])+
\mu(\alpha(\mathcal{R}(y)),[x,\mathcal{R}(z)])+\lambda\mu(\alpha(\mathcal{R}(y)),[x,z])\\&+\mu(\alpha(y),[\mathcal{R}(x),\mathcal{R}(z)])+
\lambda\mu(\alpha(y),[\mathcal{R}(x),z])+\lambda\mu(\alpha(y),[x,\mathcal{R}(z)])+\lambda^2\mu(\alpha(y),[x,z]).
\end{align*}
Finally, we find
\begin{align*}
 [\alpha(x),\mu_{\mathcal{R}}(y,z)]_{\mathcal{R}}-\mu_{\mathcal{R}}([x,y]_{\mathcal{R}},\alpha(z))&-
 (-1)^{|x||y|}\mu_{\mathcal{R}}(\alpha(y),[x,z]_{\mathcal{R}})=\Big([\alpha({\mathcal{R}}(x)),\mu(\mathcal{R}(y),z)]-
 [\mathcal{R}(x),\mathcal{R}(y)],\alpha(z))\\&-(1)^{|x||y|}\mu(\alpha(\mathcal{R}(y)),[\mathcal{R}(x),z])\Big)+
 \Big([\alpha({\mathcal{R}}(x)),\mu(y,\mathcal{R}(z))]-\mu([\mathcal{R}(x),y],\alpha(\mathcal{R}(z)))\\&-
 (-1)^{|x||y|}\mu(\alpha(y),[\mathcal{R}(x),\mathcal{R}(z)])\Big)+\Big([\alpha(x),\mu(\mathcal{R}(y),\mathcal{R}(z))]-
 \mu([x,\mathcal{R}(y)],\alpha(\mathcal{R}(z)))\\&-(-1)^{|x||y|}\mu(\alpha(\mathcal{R}(y)),[x,\mathcal{R}(z)])\Big)+
 \lambda\Big([\alpha(\mathcal{R}(x)),\mu(y,z)]-\mu([\mathcal{R}(x),y],\alpha(z))\\&-
 (-1)^{|x||y|}\mu(\alpha(y),[\mathcal{R}(x),z])\Big)+\lambda\Big([\alpha(x),\mu(\mathcal{R}(y),z)]-
 \mu([x,\mathcal{R}(y)],\alpha(z))\\&-(-1)^{|x||y|}\mu(\alpha(\mathcal{R}(y)),[x,z])\Big)+
 \lambda\Big([\alpha(x),\mu(y,\mathcal{R}(z))]-\mu([x,y],\alpha(\mathcal{R}(z)))\\&-
 (-1)^{|x||y|}\mu(\alpha(y),[x,\mathcal{R}(z)])\Big)+\lambda^2\Big([\alpha(x),\mu(y,z)]-\mu([x,y],\alpha(z))\\&-
 (-1)^{|x||y|}\mu(\alpha(y),[x,z])\Big)=0.
\end{align*}
This from the fact that $(\mathcal{A},[\cdot,\cdot],\mu,\alpha)$ is a Hom-Poisson superalgebra and we conclude that
 $(\mathcal{A},[\cdot,\cdot]_{\mathcal{R}},\mu_{\mathcal{R}},\alpha)$ is a Hom-Poisson superalgebra.
\end{proof}
 \subsection{Representation of Hom-Poisson Superalgebras. }

\begin{df}
A representation of commutative Hom-associative superalgebra $(\mathcal{A},\mu,\alpha)$ is a triple $(V,\eta,\alpha_V)$ consisting of a
$\mathbb{Z}_2$-graded vector space $V=V_{\overline{0}}\oplus V_{\overline{1}}$, an even linear maps $\eta:\mathcal{A}\to End(V)$
 and an even linear map $\alpha_V:V\to V$ satisfying for all $x,y\in\mathcal{H}(\mathcal{A})$:
\begin{eqnarray}
   \eta(\mu(x,y))\alpha_V&=&\eta(\alpha(x))\eta(y),\label{repres-Hom-ass1}\\
   \eta(\alpha(x))\eta(y)&=&(-1)^{|x||y|}\eta(\alpha(y))\eta(x)\label{repres-Hom-ass2}.
\end{eqnarray}
\end{df}
\begin{prop}\label{sumdirecthomassoc}
Let $(\mathcal{A},\mu,\alpha)$ be a commutative Hom-associative superalgebra.
Then $(V,\eta,\alpha_V)$ is a representation of $(\mathcal{A},\mu,\alpha)$ if and only if $(\mathcal{A}\oplus V,\mu_{\mathcal{A}\oplus V},\alpha+\alpha_V)$
is a commutative Hom-associative superalgebra, where $\mu_{\mathcal{A}\oplus V}$ and $(\alpha+\alpha_V)$
are defined for all $x,y\in\mathcal{H}(\mathcal{A}),\;u,v \in\mathcal{H}(V)$ by
\begin{eqnarray}
    \mu_{\mathcal{A}\oplus V}(x+u,y+v)&=&\mu(x,y)+\eta(x)v+(-1)^{|y||u|}\eta(y)u,\label{Hom-ass-direct-sum1}\\
    (\alpha+\alpha_V)(x+u)&=&\alpha(x)+\alpha_V(u)\label{Hom-ass-direct-sum2}.
\end{eqnarray}
This Hom-associative superalgebra is called semi-direct product of $(\mathcal{A},\mu,\alpha)$ and $(V,\eta,\alpha_V)$ and denoted by
 $\mathcal{A}\ltimes^{\alpha}_{\eta,\alpha_V}V$ or simply $\mathcal{A}\ltimes V$.
\end{prop}
\begin{proof}
Let $x,y,z\in\mathcal{H}(\mathcal{A}),u,v,w\in\mathcal{H}(V)$, we have
\begin{align*}
\mu_{\mathcal{A}\oplus V}(\mu_{\mathcal{A}\oplus V}(x+u,y+v),(\alpha+\alpha_V)(z+w))&=
\mu_{\mathcal{A}\oplus V}(\mu(x,y)+\eta(x)v+(-1)^{|y||u|}\eta(y)u,\alpha(z)+\alpha_V(w))\\&=
\mu(\mu(x,y),\alpha(z))+\eta(\mu(x,y))\alpha_V(w)+(-1)^{|z|(|x|+|v|)}\eta(\alpha(z))\eta(x)v\\&+(-1)^{|u|(|y|+|z|)+|y||z|}\eta(\alpha(z))\eta(y)u,
\end{align*}
and
\begin{align*}
 \mu_{\mathcal{A}\oplus V}((\alpha+\alpha_V)(x+u),\mu_{\mathcal{A}\oplus V}(y+v,z+w))&=
 \mu_{\mathcal{A}\oplus V}(\alpha(x)+\alpha_V(u),\mu(y,z)+\eta(y)w+(-1)^{|z||v|}\eta(z)v)\\&=
 \mu(\alpha(x),\mu(y,z))+\eta(\alpha(x))\eta(y)w+(-1)^{|z||v|}\eta(\alpha(x))\eta(z)v\\& +(-1)^{|u|(|y|+|z|)}\eta(\mu(y,z))\alpha_V(u).
\end{align*}
 Therefore, the Hom-associator on $\mathcal{A}\oplus V$ is  given by
\begin{align*}
ass_{\mu_{\mathcal{A}\oplus V}}(x+u,y+v,z+w)&= \mu_{\mathcal{A}\oplus V}(\mu_{\mathcal{A}\oplus V}(x+u,y+v),(\alpha+\alpha_V)(z+w))\\&-
 \mu_{\mathcal{A}\oplus V}((\alpha+\alpha_V)(x+u),\mu_{\mathcal{A}\oplus V}(y+v,z+w))\\&=
 ass_\mu(x,y,z)+\Big(\eta(\mu(x,y))\alpha_V-\eta(\alpha(x))\eta(y)\Big)(w)\\&-
 (-1)^{|z||v|}\Big(\eta(\alpha(x))\eta(z)-(-1)^{|x||z|}\eta(\alpha(z))\eta(x)\Big)(v)\\&-
 (-1)^{|u|(|y|+|z|)}\Big(\eta(\mu(y,z))\alpha_V-(-1)^{|y||z|}\eta(\alpha(z))\eta(y)\Big)(u).
\end{align*}

Since $(\mathcal{A},\mu,\alpha)$ is a commutative Hom-associative superalgebra, then $ass_\mu(x,y,z)=0$, therefore
$(\mathcal{A}\oplus V,\mu_{\mathcal{A}\oplus V},\alpha+\alpha_V)$ is a commutative Hom-associative superalgebra if and only if
\begin{align*}
ass_{\mu_{\mathcal{A}\oplus V}}(x+u,y+v,z+w)&=\Big(\eta(\mu(x,y))\alpha_V-\eta(\alpha(x))\eta(y)\Big)(w)-
(-1)^{|z||v|}\Big(\eta(\alpha(x))\eta(z)-(-1)^{|x||z|}\eta(\alpha(z))\eta(x)\Big)(v)\\&-
(-1)^{|u|(|y|+|z|)}\Big(\eta(\mu(y,z))\alpha_V-(-1)^{|y||z|}\eta(\alpha(z))\eta(y)\Big)(u)\\&=0,
\end{align*}
for all $x,y,z\in\mathcal{H}(\mathcal{A}),u,v,w\in\mathcal{H}(V)$.\\
The case $u=v=0$ gives $\eta(\mu(x,y))\alpha_V-\eta(\alpha(x))\eta(y)=0$, and the case $u=w=0$ gives
$\eta(\alpha(x))\eta(y)-(-1)^{|x||y|}\eta(\alpha(y))\eta(x)=0$ for all $x,y\in\mathcal{H}(\mathcal{A})$ which shows the proposition.
\end{proof}
\begin{df}
A representation of Hom-Lie superalgebras $(\mathcal{A},[\cdot,\cdot],\alpha)$ is a triple $(V,\rho,\alpha_V)$ where $V$ is a
$\mathbb{Z}_2$-graded vector space, $\alpha_V\in gl(V)$ and $\rho:\mathcal{A}\to End(V)$
are two even linear maps such that the following equality holds for all $x,y\in\mathcal{H}(\mathcal{A})$:
\begin{eqnarray}
\rho([x,y])\alpha_V=\rho(\alpha(x))\rho(y)-(-1)^{|x||y|}\rho(\alpha(y))\rho(x).\label{repres-Hom-Lie2}
\end{eqnarray}
\end{df}

\begin{prop}\label{sumdirecthomLie}
Let $(\mathcal{A}, [\cdot,\cdot ],\alpha)$ be a Hom-Lie superalgebra and $(V,\alpha_V)$ be a Hom-module. Let
$\rho:\mathcal{A} \to gl(V )$ be an even linear map. The triple $(V,\rho,\alpha_V)$ is a representation of
$(\mathcal{A}, [\cdot,\cdot ],\alpha)$ if and only if the direct sum  $\mathcal{A}\oplus V $ of $\mathbb{Z}_2$-graded vector spaces $\mathcal{A}$ and $V$,
turns into a Hom-Lie superalgebra by defining the linear map \eqref{Hom-ass-direct-sum2} and the folowing multiplication
\begin{equation}
    [x+u,y+v]_{\mathcal{A}\oplus V}=[x,y]+\rho(x)v-(-1)^{|y||u|}\rho(y)u,\label{Hom-Lie-direct-sum}.
\end{equation}
\end{prop}
\begin{proof}
It is easy to prove this result by using the same technique used in the proof of Proposition \ref{sumdirecthomassoc}.
\end{proof}
\begin{df}
A representation of a Hom-Poisson superalgebra $(\mathcal{A}, [\cdot,\cdot ],\mu,\alpha)$ is a quadruple
$(V,\rho,\eta,\alpha_V)$ consisting of a $\mathbb{Z}_2$-graded vector space $V$, $\alpha_V\in gl(V)$ and $\rho,\eta:\mathcal{A}\to End(V)$
are three even linear maps such that:
\begin{enumerate}
   \item $(V,\eta,\alpha_V)$ is a representation of the Hom-associative commutative superalgebras $(\mathcal{A},\mu,\alpha)$.
   \item $(V,\rho,\alpha_V)$ is a representation of the Hom-Lie superalgebras $(\mathcal{A}, [\cdot,\cdot ],\alpha)$.
   \item The following equalities holds for all $x,y\in\mathcal{H}(\mathcal{A})$:
   \begin{eqnarray}
     \eta([x,y])\alpha_V&=\rho(\alpha(x))\eta(y)-(-1)^{|x||y|}\eta(\alpha(y))\rho(x),\label{repr-Hom-Poisson1}\\
     \rho(\mu(x,y))\alpha_V&=\eta(\alpha(x))\rho(y)+-(-1)^{|x||y|}\eta(\alpha(y))\rho(x).\label{repr-Hom-Poisson2}
   \end{eqnarray}
\end{enumerate}
\end{df}

\begin{thm}
Let $(\mathcal{A}, [\cdot,\cdot ],\mu,\alpha)$ be a Hom-Poisson superalgebra and $(V,\alpha_V)$ be a Hom module.
Let $\rho,\eta:\mathcal{A}\to End(V)$ be two even linear maps. Then, $(V,\rho,\eta,\alpha_V)$ is a representation of
$(\mathcal{A}, [\cdot,\cdot ],\mu,\alpha)$ if and only if the direct sum $\mathcal{A}\oplus V $
turns into a Hom-Poisson superalgebra by defining the  linear map defined by eq.\eqref{Hom-ass-direct-sum2} and
two multiplication defined by eqs. \eqref{Hom-ass-direct-sum1}, \eqref{Hom-Lie-direct-sum}.

\end{thm}
\begin{proof}
By Propositions \ref{sumdirecthomassoc} and \ref{sumdirecthomLie} we show that
$({\mathcal{A}\oplus V},\mu_{\mathcal{A}\oplus V},(\alpha+\alpha_V))$ is a Hom-associative superalgebra
and $(\mathcal{A},[\cdot,\cdot]_{\mathcal{A}\oplus V},\alpha_{\mathcal{A}\oplus V})$ is a Hom-Lie superalgebra.\\

Now, it remains to show that the Hom-Leibniz identity is satisfied on $\mathcal{A}\oplus V$.
Let $x,y,z\in\mathcal{H}(\mathcal{A})$ and $u,v,w\in\mathcal{H}(V)$, we have
\begin{align*}
  \bullet)\; [(\alpha+\alpha_V) (x+u),\mu_{\mathcal{A}\oplus V}(y+v,z+w)]_{\mathcal{A}\oplus V}&=
  [\alpha(x)+\alpha_V(u),\mu(y,z)+\eta(y)w+\eta(z)v,]_{\mathcal{A}\oplus V}\\&=
  [\alpha(x),\mu(y,z)]+\rho(\alpha(x))\eta(y)w+\rho(\alpha(x))\eta(z)v-\rho(\mu(y,z))\alpha_V(u).
\end{align*}
\begin{align*}
   \bullet)\; \mu_{\mathcal{A}\oplus V}\big([x+u,y+v]_{\mathcal{A}\oplus V},(\alpha+\alpha_V)(z+w)\big)&=
   \mu_{\mathcal{A}\oplus V}([x,y]+\rho(x)v-\rho(y)u,\alpha(z)+\alpha_V(w))\\&=\mu([x,y],\alpha(z))+\eta([x,y])\alpha_V(w))
   +\eta(\alpha(z))\rho(x)v-\eta(\alpha(z))\rho(y)u.
\end{align*}
\begin{align*}
\bullet)\; \mu_{\mathcal{A}\oplus V}\big(\alpha_{\mathcal{A}\oplus V}(y+v),[x+u,z+w]_{\mathcal{A}\oplus V}\big)&=
\mu_{\mathcal{A}\oplus V}(\alpha(y)+\alpha_V(v),[x,z]+\rho(x)w-\rho(z)u)\\&=\mu(\alpha(y),\mu(x,z))+
\eta(\alpha(y))\rho(x)w-\eta(\alpha(y))\rho(z)u)+\eta([x,z])\alpha_V(v).
\end{align*}
By a direct computation, the Hom-Leibniz identity is satisfied on $\mathcal{A}\oplus V$ if and only if the eqs.
 \eqref{repr-Hom-Poisson1} and \eqref{repr-Hom-Poisson2} are satisfied which implies that $(V,\rho,\eta,\alpha_V)$
 is a representation of $(\mathcal{A}, [\cdot,\cdot ],\mu,\alpha)$.
\end{proof}

 \section{Some construction of $n$-ary Hom-Nambu Poisson Superalgebras. }

 \subsection{$n$-ary Hom-Nambu Poisson Superalgebras. }
\begin{df}
An \emph{$n$-ary Hom-Nambu} superalgebra $(\mathcal{N}, [\cdot,\dots,\cdot],  \widetilde{\alpha} )$ is a triple consisting of a $\mathbb{Z}_2$-graded linear space
   $\mathcal{N}=\mathcal{N}_0\oplus\mathcal{N}_1$, an even
$n$-linear map $[\cdot ,\dots, \cdot ]:  \mathcal{N}^{ n}\to \mathcal{N}$ ( i.e. $[\mathcal{N}_{k_1} ,\dots,\mathcal{N}_{k_n}]
 \subset \mathcal{N}_{k_1+\dots +k_n}$) and a family
$\widetilde{\alpha}=(\alpha_i)_{1\leq i\leq n-1}$ of even linear maps $ \alpha_i:\ \ \mathcal{N}\to \mathcal{N}$, satisfying
  \begin{eqnarray}\label{NambuIdentity}
&& \forall (x_1,\dots, x_{n-1})\in \mathcal{H}(\mathcal{N})^{n-1}, \,\, (y_1,\dots,  y_n)\in \mathcal{H}(\mathcal{N})^{ n}: \nonumber \\
&& \big[\alpha_1(x_1),\dots.,\alpha_{n-1}(x_{n-1}),[y_1,\dots.,y_{n}]\big]= \\ \nonumber
&& \sum_{i=1}^{n}(-1)^{|X||Y|^{i-1}}\big[\alpha_1(y_1),\dots.,\alpha_{i-1}(y_{i-1}),[x_1,\dots.,x_{n-1},y_i], 
\alpha_i(y_{i+1}),\dots,\alpha_{n-1}(y_n)\big],
  \end{eqnarray}
where
$|X|=\displaystyle\sum_{k=1}^{n-1}|x_k|$ and $|Y|^{i-1}=\displaystyle\sum_{k=1}^{i-1}|y_k|.$

The identity \eqref{NambuIdentity} is called \emph{super-Hom-Nambu identity}.
  \end{df}
Let $\widetilde{\alpha}:\mathcal{N}^{n-1}\to\mathcal{N}^{n-1}$ be even linear maps defined for all $X=(x_1,\ldots,x_{n-1})\in \mathcal{N}^{n-1}$ by
$\widetilde{\alpha}(X)=(\alpha_1(x_1),\ldots,\alpha_{n-1}(x_{n-1}))\in\mathcal{N}^{n-1}$.
For all $X=(x_1,\ldots,x_{n-1})\in \mathcal{N}^{n-1}$,  the map $\text{ad}_X:\mathcal{N}\to\mathcal{N}$ defined by
\begin{equation}\label{adjointMapNaire}
\text{ad}_X(y)=[x_{1},\dots,x_{n-1},y],\quad  \forall y\in \mathcal{N},
\end{equation}
is called adjoint map. Then the super-Hom-Nambu identity \eqref{NambuIdentity} may be written in terms of adjoint map as
\begin{eqnarray*}
\text{ad}_{\widetilde{\alpha} (X)}( [y_1,\dots,y_n])=
\sum_{i=1}^{n-1}(-1)^{|X||Y|^{i-1}} \left[\alpha_1(y_1),\dots,\alpha_{i-1}(y_{i-1}),
\text{ad}_X(y_{i}), \alpha_{i}(y_{i+1}) \dots,\alpha_{n-1}(y_n)\right].
\end{eqnarray*}

\begin{df}
  An $n$-ary Hom-Nambu superalgebra $(\mathcal{N}, [\cdot,\dots,\cdot],  \widetilde{\alpha} )$ is called $n$-Hom-Lie superalgebra if the bracket
  $[\cdot,\dots,\cdot]$ is super-skewsym\-metric that is
  \begin{align}\label{SuperSkewSym}
& [x_1,\dots,x_i,x_{i+1},\dots,x_n]=-(-1)^{|x_i||x_{i+1}|}[x_1,\dots,x_{i+1},x_i,\dots,x_n],\forall\; 1\leq i\leq n-1.
  \end{align}
  It is equivalent to
 \begin{align}\label{SuperSkewSym1}
&\left[x_1,\dots,x_i,\dots,x_j,\dots,x_n \right]=-(-1)^{|X|^{j-1}_{i+1}(|x_i|+|x_j|)+|x_i||x_j|}  \left[x_1,\dots,x_j,\dots,x_i,\dots,x_n\right],\;\forall\;1\leq i<j\leq n,
 \end{align}
 where $x_1,\dots ,x_n\in \mathcal{H}(\mathcal{N})$ and $|X|_{i}^{j}=\displaystyle\sum_{k=i}^{j}|x_k|.$
\end{df}
\begin{rem}
When the maps $(\alpha_i)_{1\leq i\leq n-1}$ are all identity maps, one recovers the classical $n$-ary Nambu superalgebras.
\end{rem}

Let $(\mathcal{N},[\cdot,\dots,\cdot],\widetilde{\alpha})$ and
$(\mathcal{N}',[\cdot,\dots,\cdot]',\widetilde{\alpha}')$ be two $n$-ary Hom-Nambu
superalgebras  where $\widetilde{\alpha}=(\alpha_{i})_{1\leq i\leq n-1}$ and
$\widetilde{\alpha}'=(\alpha'_{i})_{1\leq i\leq n-1}$. A linear map $f:\mathcal{N}\to \mathcal{N}'$ is an
$n$-ary Hom-Nambu superalgebras \emph{morphism}  if it satisfies
\begin{eqnarray*}&f([x_{1},\dots,x_n])=[f(x_{1}),\dots,f(x_n)]',\\
&f \circ \alpha_i=\alpha'_i\circ f, \quad \forall i=1,\dots,n-1.
\end{eqnarray*}

In the sequel we deal sometimes with a particular class of $n$-ary Hom-Nambu superalgebras which we call $n$-ary multiplicative Hom-Nambu  superalgebras.

\begin{df}
A multiplicative $n$-ary Hom-Nambu superalgebra 
(resp.  multiplicative $n$-Hom-Lie superalgebra) is an $n$-ary Hom-Nambu superalgebra  (resp. $n$-Hom-Lie superalgebra)
$(\mathcal{N}, [\cdot,\dots,\cdot],  \widetilde{ \alpha})$, with  $\widetilde{\alpha}=(\alpha_i)_{1\leq i\leq n-1}$
where  $\alpha_1= \dots =\alpha_{n-1}=\alpha$  and satisfying
\begin{equation}
\alpha([x_1,\dots,x_n])=[\alpha(x_1),\dots,\alpha(x_n)],\ \  \forall x_1,\dots,x_n\in \mathcal{N}.
\end{equation}
For simplicity, denote the $n$-ary multiplicative Hom-Nambu superalgebra as $(\mathcal{N}, [\cdot,\dots,\cdot],  \alpha)$
 where $\alpha :\mathcal{N}\to \mathcal{N}$ is an even linear map. Also by misuse of language an element  $X\in \mathcal{N}^n$
  refers to  $X=(x_1,\dots,x_{n})$ and  $\alpha(X)$ denotes
$(\alpha (x_1),\dots,\alpha (x_n))$.
\end{df}
\begin{df}
A multiplicative $n$-ary Hom-Nambu superalgebra $(\mathcal{N}, [\cdot ,\dots, \cdot],  \alpha)$ is called regular if $\alpha$ is bijective.
\end{df}
\begin{df}
An $n$-ary Hom-Nambu Poisson superalgebras is a quadruple $(\mathcal{A},[\cdot,\cdots,\cdot],\mu,\alpha)$ consisting of a $\mathbb{Z}_2$-graded linear space
$\mathcal{A}=\mathcal{A}_0\oplus\mathcal{A}_1$, an even $n$-linear map $[.,\cdots,.]:\mathcal{A}^n\rightarrow\mathcal{A}$, an
even bilinear map $\mu:\mathcal{A}\times\mathcal{A}\rightarrow\mathcal{A}$ and an even linear map $\alpha:\mathcal{A}\rightarrow\mathcal{A}$
such that:
\begin{enumerate}
\item $(\mathcal{A},\mu,\alpha)$ is a Hom-associative commutative superalgebra.
\item $(\mathcal{A},[\cdot,\cdots,\cdot],\alpha)$ is an $n$-Hom-Lie superalgebra.
\item The following identity holds for all $x_1,\cdots,x_{n-1},y,z\in\mathcal{H}(\mathcal{A})$:
\begin{equation}\label{n-aryHomLeibniz}
[\alpha(x_1),\cdots,\alpha(x_{n-1}),\mu(y,z)]=(-1)^{|y||X|}\mu\big(\alpha(y),[x_1,\cdots,x_{n-1},z]\big)+
\mu\big([x_1,\cdots,x_{n-1},y],\alpha(z)\big),
\end{equation}
where $|X|=\displaystyle\sum_{i=0}^{n-1}|x_i|$.
\end{enumerate}
The identity $\eqref{n-aryHomLeibniz}$ is called Hom-Leibniz identity.
\end{df}
\begin{rem}
The commutativity of $\mu$ gives that the identity $\eqref{n-aryHomLeibniz}$ is equivalent to
\begin{equation}\label{equiv-n-aryHomLeibniz}
[\alpha(x_1),\cdots,\alpha(x_{n-1}),\mu(y,z)]=(-1)^{|y||z|}\mu\big([x_1,\cdots,x_{n-1},z],\alpha(y)\big)+
\mu\big([x_1,\cdots,x_{n-1},y],\alpha(z)\big),
\end{equation}
\end{rem}

 \subsection{n-ary Hom-Nambu Poisson Superalgebras induced by Hom-Poisson Superalgebras }.\\

In \cite{Mabrouk-Ncib-Sergei}, the authors introduced a construction of $n$-Hom-Lie superalgebra by a Hom-Lie superalgebra
 in which they gave an $n$-ary product defined by
 \begin{equation}\label{nProduct}
[x_1,\dots,x_n]_\phi=\sum_{i<j}^{n}(-1)^{i+j+1}(-1)^{\gamma^X_{ij}}\phi(x_1,\dots,\hat{x_i},\dots,\hat{x_j},\dots,x_n)[x_i,x_j],
\end{equation}
for all $x_1,\dots,x_n\in\mathcal{H}(\mathcal{A})$, where $\gamma^X_{ij}=|X|^n_{j+1}(|x_i|+|x_j|)+|x_i||X|_{j-1}^{i+1}$ and $\phi\in \wedge^{n-2}\mathcal{A}^*$ be an even $(n-1)$-cochain
. The authors showed
that this product defines a structure of a multiplicative n-Hom-Lie superalgebra under some properties as in the following theorem.
\begin{thm}\cite{Mabrouk-Ncib-Sergei}\label{HomLieSupAlToNaryHomLieSupAlg}
  Let $(\mathcal{A},[\cdot,\cdot],\alpha)$ be a  multiplicative Hom-Lie superalgebra, $\mathcal{A}^*$  its dual and $\phi$ be
an even cochain of degree $n-2$, i.e. $\phi\in\wedge^{n-2}\mathcal{A}^*$. The linear space   $\mathcal{A}$
equipped with the n-ary product \eqref{nProduct} and the even linear map $\alpha$ is a multiplicative n-Hom-Lie superalgebra if and only if
\begin{align}\label{NHomLieProduct}
 & \phi\wedge\delta\phi_X=0,\ \forall X\in \wedge^{n-3}\mathcal{H}(\mathcal{A}),\\
  & \phi\circ(\alpha\otimes Id\otimes\dots\otimes Id)=\phi.\label{NHomLieProduct1}
\end{align}
\end{thm}
In this section we extend this construction to the $n$-ary-Hom-poisson superalgebras.
\begin{thm}
Let $(\mathcal{A},[\cdot,\cdot],\mu,\alpha)$ be a multiplicative Hom-Poisson superalgebra and $\phi\in\wedge^{n-2}\mathcal{A}^*$ satisfying the conditions
 \eqref{NHomLieProduct} and \eqref{NHomLieProduct1}. Then, the quadruple $(\mathcal{A},[\cdot,\cdots,\cdot]_\phi,\mu,\alpha)$
 is a multiplicative $n$-ary Hom-Poisson superalgebra if and only if
\begin{align}
  & \phi(X,\mu(y,z))\alpha(t)=\phi(X,y)\mu(\alpha(z),t)+(-1)^{|y||z|}\phi(X,z)\mu(\alpha(y),t),\label{NHomPoisson2}
\end{align}
for all $X\in\wedge^{n-3}\mathcal{H}(\mathcal{A}),\;y,z,t\in\mathcal{H}(\mathcal{A})$, which is called the induce $n$-ary Hom-Poisson superalgebra of $(\mathcal{A},[\cdot,\cdot],\mu,\alpha)$.
\end{thm}

\begin{proof}
By the Theorem \ref{HomLieSupAlToNaryHomLieSupAlg},  $(\mathcal{A},[\cdot,\cdots,\cdot]_\phi,\alpha)$
is a multiplicative $n$-Hom-Lie superalgebra. It remains to show the Hom-Leibniz identity.\\
Let $X=(x_1,\cdots,x_{n-1})\in\wedge^{n-1}\mathcal{H}(\mathcal{A}),y,z\in\mathcal{H}(\mathcal{A})$, we have
{\small\begin{align*}
  [\alpha(x_1),\cdots,\alpha(x_{n-1}),\mu(y,z)]_\phi&=\displaystyle\sum_{i<j\leq n-1}(-1)^{i+j+1}
  (-1)^{\gamma^{X^{ij}}_{yz}}\phi(\alpha(x_1),\cdots,\hat{\alpha(x_i)},\cdots,\hat{\alpha(x_j)},\cdots,\alpha(x_{n-1}),\mu(y,z))[\alpha(x_i),\alpha(x_j)] \\
 &+\displaystyle\sum_{i=1}^{n-1}(-1)^{i+n+1}(-1)^{|x_i||X|_{i+1}}
 \phi(\alpha(x_1),\cdots,\hat{\alpha(x_i)},\cdots,\alpha(x_{n-1}),\hat{\mu(y,z)})[\alpha(x_i),\mu(y,z)],
\end{align*}}
where $$\gamma^{X^{ij}}_{yz}=(|x_i|+|x_j|)(|x_{j+1}|+\cdots+|x_{n-1}|+|y|+|z|)+|x_i|(|x_{i+1}|+\cdots+|x_{j-1}|)
~~and~~|X|_{k}=|x_{k}|+\cdots+|x_{n-1}|,\forall1\leq k\leq n-1.$$
On the other hand, we get
\begin{align*}
   \mu([x_1,\cdots,x_{n-1},z]_\phi,\alpha(y))&=\displaystyle\sum_{i<j\leq n-1}(-1)^{i+j+1}(-1)^{\gamma^{X^{ij}}_{z}}
   \phi(x_1,\cdots,\hat{x_i},\cdots,\hat{x_j},\cdots,x_{n-1},z)\mu([x_i,x_j],\alpha(y))\\&+\displaystyle\sum_{i=1}^{n-1}
   (-1)^{i+n+1}(-1)^{|x_i||X|_{i+1}}\phi(x_1,\cdots,\hat{x_i},\cdots,x_{n-1},\hat{z})\mu([x_i,z],\alpha(y)),
\end{align*}
where $$\gamma^{X^{ij}}_z=(|x_i|+|x_j|)(|x_{j+1}|+\cdots+|x_{n-1}|+|z|)+|x_i|(|x_{i+1}|+\cdots+|x_{j-1}|)$$ and
\begin{align*}
    \mu([x_1,\cdots,x_{n-1},y]_\phi,\alpha(z))&=\displaystyle\sum_{i<j\leq n-1}(-1)^{i+j+1}(-1)^{\gamma^{X^{ij}}_{y}}
    \phi(x_1,\cdots,\hat{x_i},\cdots,\hat{x_j},\cdots,x_{n-1},y)\mu([x_i,x_j],\alpha(z))\\&+\displaystyle\sum_{i=1}^{n-1}
    (-1)^{i+n+1}(-1)^{|x_i||X|_{i+1}}\phi(x_1,\cdots,\hat{x_i},\cdots,x_{n-1},\hat{y})\mu([x_i,y],\alpha(z)).
\end{align*}
Since $(\mathcal{A},[\cdot,\cdot],\mu,\alpha)$ is a multiplicative Hom-Lie superalgebra and by the identity
\eqref{equiv-Homleibnizident}, we have
$$[\alpha(x_i),\mu(y,z)]=\mu([x_i,y],\alpha(z))+(-1)^{|y||z|}\mu([x_i,z],\alpha(y)),$$
which gives
\begin{align*}
[\alpha(x_1),\cdots,\alpha(x_{n-1}),\mu(y,z)]_\phi&-(-1)^{|y||z|}\mu([x_1,\cdots,x_{n-1},z]_\phi,\alpha(y))-
\mu([x_1,\cdots,x_{n-1},y]_\phi,\alpha(z))\\&=\displaystyle\sum_{i<j\leq n-1}(-1)^{i+j+1}
(-1)^{\gamma^{X^{ij}}_{yz}}\phi(\alpha(x_1),\cdots,\hat{\alpha(x_i)},\cdots,\hat{\alpha(x_j)},\cdots,\alpha(x_{n-1}),\mu(y,z))[\alpha(x_i),\alpha(x_j)]
\\&-(-1)^{|y||z|}\displaystyle\sum_{i<j\leq n-1}(-1)^{i+j+1}(-1)^{\gamma^{X^{ij}}_{z}}
\phi(x_1,\cdots,\hat{x_i},\cdots,\hat{x_j},\cdots,x_{n-1},z)\mu([x_i,x_j],\alpha(y))\\&-
\displaystyle\sum_{i<j\leq n-1}(-1)^{i+j+1}(-1)^{\gamma^{X^{ij}}_{y}}\phi(x_1,\cdots,\hat{x_i},\cdots,\hat{x_j},\cdots,x_{n-1},y)\mu([x_i,x_j],\alpha(z))
\\&=\displaystyle\sum_{i<j\leq n-1}(-1)^{i+j+1}(-1)^{\gamma^{X^{ij}}_{yz}}
\phi(x_1,\cdots,\hat{\alpha(x_i)},\cdots,\hat{\alpha(x_j)},\cdots,x_{n-1},\mu(y,z))\alpha([x_i,x_j])\\&-
\displaystyle\sum_{i<j\leq n-1}(-1)^{i+j+1}(-1)^{\gamma^{X^{ij}}_{yz}}(-1)^{|y||z|}
\phi(x_1,\cdots,\hat{x_i},\cdots,\hat{x_j},\cdots,x_{n-1},z)\mu(\alpha(y),[x_i,x_j])\\&-
\displaystyle\sum_{i<j\leq n-1}(-1)^{i+j+1}(-1)^{\gamma^{X^{ij}}_{yz}}
\phi(x_1,\cdots,\hat{x_i},\cdots,\hat{x_j},\cdots,x_{n-1},y)\mu(\alpha(z),[x_i,x_j])\\&=
\displaystyle\sum_{i<j\leq n-1}(-1)^{i+j+1}(-1)^{\gamma^{X^{ij}}_{yz}}
\Big(\phi(x_1,\cdots,\hat{\alpha(x_i)},\cdots,\hat{\alpha(x_j)},\cdots,x_{n-1},\mu(y,z))\alpha([x_i,x_j])\\&-
(-1)^{|y||z|}\phi(x_1,\cdots,\hat{x_i},\cdots,\hat{x_j},\cdots,x_{n-1},z)\mu(\alpha(y),[x_i,x_j])\\&-
\phi(x_1,\cdots,\hat{x_i},\cdots,\hat{x_j},\cdots,x_{n-1},y)\mu(\alpha(z),[x_i,x_j]) \Big).
\end{align*}
Suppose that $X=(x_1,\cdots,\hat{x_i},\cdots,\hat{x_j},\cdots,x_{n-1}),\;t=[x_i,x_j]$ and we apply the identity \eqref{NHomPoisson2}, we find
$$[\alpha(x_1),\cdots,\alpha(x_{n-1}),\mu(y,z)]_\phi-(-1)^{|y||z|}\mu([x_1,\cdots,x_{n-1},z]_\phi,\alpha(y))-\mu([x_1,\cdots,x_{n-1},y]_\phi,\alpha(z))=0,$$
therefore the Hom-Leibniz identity is verified and the quadruple $(\mathcal{A},[\cdot,\cdots,\cdot]_\phi,\mu,\alpha)$
has a structure of multiplicative $n$-ary Hom-Poisson superalgebra.
\end{proof}
We are interesting to see the opposite problem, i.e to construct a Hom-Poisson superalgebra from an $n$-ary Hom-Poisson superalgebra.
 Let $(\mathcal{A},[\cdot,\cdots,\cdot],\mu,\alpha)$ be an $n$-ary Hom-Poisson superagebra.
  We fix $A=(a_1,\cdots, a_{n-2})$ in $\in\mathcal{H}(\mathcal{A})^{n-2}$ and we define on $\mathcal{A}$
  the following bilinear map, for all $x,y\in\mathcal{H}(\mathcal{A})$
\begin{equation}\label{hom-poisson-induce}
[x,y]_A=[x,y,a_1,\cdots,a_{n-2}].
\end{equation}
\begin{thm}
By the above notations, the quadruple $(\mathcal{A},[\cdot,\cdot]_A,\mu,\alpha)$ define a Hom-Poisson superalgebra if and only if
\begin{align}\label{induiceHomLieProduct1}
 & a_i\in\ker(\alpha-Id),\ \forall i\in \{1,\cdots,n-2\},\\
  & |A|=\displaystyle\sum_{i=1}^{n-2}|a_i|=0.\label{induceHomLieProduct2}
\end{align}
\end{thm}
\begin{proof}
It's easy to see that $[\cdot,\cdot]_A$ is super-skew-symmetry.\\

Let $x,y,z\in\mathcal{H}(\mathcal{A})$. By the definition of $[\cdot,\cdot]_A$ and the equations \eqref{induiceHomLieProduct1}-\eqref{induceHomLieProduct2},
 we have
\begin{align*}
[\alpha(x),[y,z]_A]_A&=[\alpha(x),[y,z,a_1,\cdots,a_{n-2}],a_1,\cdots,a_{n-2}]\\&=
(-1)^{|A|(|A|+|y|+|z|)}[\alpha(x),a_1,\cdots,a_{n-2},[y,z,a_1,\cdots,a_{n-2}]]\\&=
[\alpha(x),\alpha(a_1),\cdots,\alpha(a_{n-2}),[y,z,a_1,\cdots,a_{n-2}]]\\&=
[[x,a_1,\cdots,a_{n-2},y],\alpha(z),\alpha(a_1),\cdots,\alpha(a_{n-2})]\\&+
(-1)^{|y|(|A|+|x|)}[\alpha(y),[x,a_1,\cdots,a_{n-2},z],\alpha(a_1),\cdots,\alpha(a_{n-2})]\\&+
\displaystyle\sum_{i=1}^{n-2}(-1)^{|A|^{i-1}(|A|+|x|)}
[\alpha(y),\alpha(z),\alpha(a_1),\cdots,\alpha(a_{i-1}),[x,a_1,\cdots,a_{n-2},a_i],\alpha(a_{i+1}),\cdots,\alpha(a_{n-2})].
\end{align*}
The last equality is since the equation \eqref{NambuIdentity} and for all $i\in\{1,\cdots,n-2\},\;[x,a_1,\cdots,a_{n-2},a_i]=0$, which gives
\begin{align*}
[\alpha(x),[y,z]_A]_A&=[[x,a_1,\cdots,a_{n-2},y],\alpha(z),\alpha(a_1),\cdots,\alpha(a_{n-2})]\\&+
(-1)^{|x||y|}[\alpha(y),[x,a_1,\cdots,a_{n-2},z],\alpha(a_1),\cdots,\alpha(a_{n-2})]\\&=
[[x,a_1,\cdots,a_{n-2},y],\alpha(z),a_1,\cdots,a_{n-2}]+(-1)^{|x||y|}[\alpha(y),[x,a_1,\cdots,a_{n-2},z],a_1,\cdots,a_{n-2}]\\&=
(-1)^{|A||y|}[[x,y,a_1,\cdots,a_{n-2}],\alpha(z),a_1,\cdots,a_{n-2}]\\&+(-1)^{|x||y|+|A||z|}
[\alpha(y),[x,z,a_1,\cdots,a_{n-2}],a_1,\cdots,a_{n-2}]\\&=[[x,y]_A,\alpha(z)]_A+(-1)^{|x||y|}[\alpha(y),[x,z]_A]_A.
\end{align*}
Therefore the Hom super-Jacobi identity is satisfied and the triple $(\mathcal{A},[\cdot,\cdot]_A,\alpha)$ is a Hom-Lie superalgebras.\\

It remains to show the Hom-Leibniz identity on $(\mathcal{A},[\cdot,\cdot]_A,\mu,\alpha)$.
By using the equation \eqref{induceHomLieProduct2} and the Hom-Leibniz identity on $(\mathcal{A},[\cdot,\cdots,\cdot],\mu,\alpha)$, we have
\begin{align*}
  [x,\mu(y,z)]_A&=[x,\mu(y,z),a_1,\cdots,a_{n-2}]\\&=[x,a_1,\cdots,a_{n-2},\mu(y,z)]\\&=
  \mu([x,a_1,\cdots,a_{n-2},y],\alpha(z))+(-1)^{|y|(|A|+|x|)}  \mu(\alpha(y),[x,a_1,\cdots,a_{n-2},z])\\&=
  \mu([x,y,a_1,\cdots,a_{n-2}],\alpha(z))+(-1)^{|x||y|}  \mu(\alpha(y),[x,z,a_1,\cdots,a_{n-2}])\\&=
  \mu([x,y]_A,\alpha(z))+(-1)^{|x||y|}\mu(\alpha(y),[x,z]_A),
\end{align*}
for all $x,y,z\in\mathcal{H}(\mathcal{A})$. So the Hom-Leibniz identity on  $(\mathcal{A},[\cdot,\cdot]_A,\mu,\alpha)$ is satisfied.\\
The theorem is proved.
\end{proof}

 \section{Representation of n-ary Hom-Nambu Poisson Superalgebras. }
\begin{df}
A representation of an $n$-Hom-Lie Superalgebra $(\mathcal{A},[\cdot,\cdots,\cdot],\alpha)$ is a triple $(V,\rho,\phi)$ consisting of a
$\mathbb{Z}_2$-graded vector space $V$, an even linear map $\rho:\Lambda^{n-1}\mathcal{A}\to gl(V)$ and an even linear map $\alpha_V:V\to V$
such that for all $x_1,\cdots,x_{n-1},y_1,\cdots,y_n\in\mathcal{H}(\mathcal{A})$, we have
\begin{align}
    &\rho(\alpha(x_1),\cdots,\alpha(x_{n-1}))\rho(y_1,\cdots,y_{n-1})-(-1)^{|X||Y|^{n-1}}
    \rho(\alpha(y_1),\cdots,\alpha(y_{n-1}))\rho(x_1,\cdots,x_{n-1})\nonumber\\&=\displaystyle\sum_{i=1}^{n-1}(-1)^{|X||Y|^{i-1}}
    \rho(\alpha(y_1),\cdots,\alpha(y_{i-1}),[x_1,\cdots,x_{n-1},y_i],\alpha(y_{i+1}),\cdots\alpha(y_{n-1}))\alpha_V;\label{repr-n-Hom-Lie1}\\&
    \rho(\alpha(x_1),\cdots,\rho(x_{n-2}),[y_1,\cdots,y_n])\alpha_V\nonumber\\&=
    \displaystyle\sum_{i=1}^{n}(-1)^{n-i}(-1)^{|X|^{n-2}(|Y|+|y_i|)+|y_i||Y|_{i+1}}\rho(\alpha(y_1),\cdots,\hat{y_i},\cdots,\alpha(y_n))\rho(x_1,\cdots,x_{n-2},y_i)\label{repr-n-Hom-Lie2}.
\end{align}
\end{df}
\begin{prop}\label{direct-sum-n-hom-lie-superalgebras}
Let $(\mathcal{A},[\cdot,\cdots,\cdot],\alpha)$ be an $n$-Hom-Lie superalgebra, $\rho:\Lambda^{n-1}\mathcal{A}\to gl(V)$  and $\alpha_V:V\to V$ are two even linear map. Then, $(\mathcal{A}\oplus V,[\cdot,\cdots,\cdot]_{\mathcal{A}\oplus V},\alpha+\alpha_V)$ is an $n$-Hom-Lie superalgebra if and only if $(V,\rho,\alpha_V)$ is a representation of $(\mathcal{A},[\cdot,\cdots,\cdot],\alpha)$, where 
\begin{equation}\label{crochet-direct-sum-n-hom-lie-superalgebras}
    [x_1+a_1,\cdots,x_n+a_n]_{\mathcal{A}\oplus V}=[x_1,\cdots,x_n]+\displaystyle\sum_{k=1}^n(-1)^{|x_k||X|_{k+1}}\rho(x_1,\cdots,\hat{x_k},\cdots,x_n)a_k,
\end{equation}
for all $x_1,\cdots,x_n\in\mathcal{H}(\mathcal{A})$ and $a_1,\cdots,a_n\in\mathcal{H}(V)$.
\end{prop}
\begin{proof}
Let $x_1,\cdots,x_{n-1},y_1,\cdots,y_n\in\mathcal{H}(\mathcal{A})$ and $a_1,\cdots,a_{n-1},b_1,\cdots,b_n\in\mathcal{H}(V)$.\\ In the one hand, we have
\begin{align*}
&\big[(\alpha+\alpha_V)(x_1+a_1),\cdots,(\alpha+\alpha_V)(x_{n-1}+a_{n-1}),[y_1+b_1,\cdots,y_n+b_n]_{\mathcal{A}\oplus V}\big]_{\mathcal{A}\oplus V}\\&=\big[\alpha(x_1)+\alpha_V(a_1),\cdots,\alpha(x_{n-1})+\alpha_V(a_{n-1}),[y_1+b_1,\cdots,y_n+b_n]_{\mathcal{A}\oplus V}\big]_{\mathcal{A}\oplus V}\\&=\big[\alpha(x_1)+\alpha_V(a_1),\cdots,\alpha(x_{n-1})+\alpha_V(a_{n-1}),[y_1,\cdots,y_n]+\displaystyle\sum_{i=1}^n(-1)^{|y_i||Y|_{i+1}}\rho(y_1,\cdots,\hat{y_i},\cdots,y_n)b_i\big]_{\mathcal{A}\oplus V}\\&=[\alpha(x_1),\cdots,\alpha(x_{n-1}),[y_1,\cdots,y_n]]+\displaystyle\sum_{j=1}^{n-1}(-1)^{|x_j|(|X|_{j+1}|+|Y|)}\rho(\alpha(x_1),\cdots,\hat{\alpha(x_j)},\cdots,\alpha(x_{n-1}),[y_1,\cdots,y_n])\alpha_V(a_j)\\&+\displaystyle\sum_{i=1}^n(-1)^{|y_i||Y|_{i+1}} \rho(\alpha(x_1),\cdots,\alpha(x_{n-1}))\rho(y_1,\cdots,\hat{y_i},\cdots,y_n)b_i.  
\end{align*}
In the other hand, we have
\begin{align*}
    &\displaystyle\sum_{k=1}^n(-1)^{|X||Y|^{k-1}}[\alpha(y_1)+\alpha_V(b_1),\cdots,[x_1+a_1,\cdots,x_{n-1}+a_{n-1},y_k+b_k]_{\mathcal{A}\oplus V},\cdots,\alpha(y_n)+\alpha_V(b_n)]_{\mathcal{A}\oplus V}\\&=\displaystyle\sum_{k=1}^n(-1)^{|X||Y|^{k-1}}\Big[\alpha(y_1)+\alpha_V(b_1),\cdots,[x_1,\cdots,x_{n-1},y_k]+\displaystyle\sum_{i=1}^{n-1}(-1)^{|x_i|(|X|_{i+1}+|y_k|)}\rho(x_1,\cdots,\hat{x_i},\cdots,x_{n-1},y_k)a_i\\&+\rho(x_1,\cdots,x_{n-1})b_k,\cdots,\alpha(y_n)+\alpha_V(b_n)\Big]_{\mathcal{A}\oplus V}\\&=\displaystyle\sum_{k=1}^n(-1)^{|X||Y|^{k-1}}\Big[\alpha(y_1),\cdots[x_1,\cdots,x_{n-1},y_k],\cdots,\alpha(y_n)\Big]\\&+\displaystyle\sum_{k=1}^n(-1)^{|X||Y|^{k-1}}\displaystyle\sum_{j=1}^{k-1}(-1)^{|y_j|(|X|+|Y|_{j+1})}\rho\big(\alpha(y_1),\cdots,\hat{\alpha(y_j)},\cdots,[x_1,\cdots,x_{n-1},y_k],\cdots,\alpha(y_n)\big)\alpha_V(b_j)\\&+\displaystyle\sum_{k=1}^n(-1)^{|X||Y|^{k-1}}\displaystyle\sum_{j=k+1}^{n}(-1)^{|y_j||Y|_{j+1}}\rho\big(\alpha(y_1),\cdots,[x_1,\cdots,x_{n-1},y_k],\cdots,\hat{\alpha(y_j)},\cdots,\alpha(y_n)\big)\alpha_V(b_j)\\&+\displaystyle\sum_{k=1}^n(-1)^{|X||Y|^{k-1}}(-1)^{(|X|+|y_k|)|Y|_{k+1}}\rho\big(\alpha(y_1),\cdots,\alpha(y_{k-1}),\alpha(y_{k+1}),\cdots,\alpha(y_n)\big)\displaystyle\sum_{i=1}^n(-1)^{|x_i|(|X|_{i+1}+|y_k|)}\rho(x_1,\cdots,\hat{x_i},\cdots,x_{n-1},y_k)a_i\\&+\displaystyle\sum_{k=1}^n(-1)^{|X||Y|^{k-1}}(-1)^{(|X|+|y_k|)|Y|_{k+1}}\rho(\alpha(y_1),\cdots,\alpha(y_{k-1}),\alpha(y_{k+1}),\cdots,\alpha(y_n))\rho(x_1,\cdots,x_{n-1})b_k.
\end{align*}
Using eq. \eqref{NambuIdentity} and conditions of representation defined by eqs. \eqref{repr-n-Hom-Lie1},  \eqref{repr-n-Hom-Lie2}, we can show that $(\mathcal{A}\oplus V,[\cdot,\cdots,\cdot]_{\mathcal{A}\oplus V},\alpha+\alpha_V)$ is an $n$-Hom-Lie superalgebra if and only if $(V,\rho,\alpha_V)$ is a representation of $(\mathcal{A},[\cdot,\cdots,\cdot],\alpha)$.
\end{proof}

 Let $(V,\rho,\eta,\alpha_V)$ be a representation of a Hom-Poisson superalgebra $(\mathcal{A},[\cdot,\cdot],\mu,\alpha)$. We define the following even linear map $\rho_\phi:\Lambda^{n-1}\mathcal{A}\to gl(V)$ by
 \begin{equation}\label{n-ary-hom-poisson-induce-rep}
     \rho_\phi(x_1,\cdots,x_{n-1})=\displaystyle\sum_{i=1}^{n-1}(-1)^{n-i-1}(-1)^{|x_i||X|_{i+1}}\phi(x_1,\cdots,\hat{x_i},\cdots,x_{n-1})\rho(x_i),
 \end{equation}
 for all $x_1,\cdots,x_{n-1}\in\mathcal{H}(\mathcal{A})$, where $\phi\in\wedge^{n-2}\mathcal{A}^*$ and $|X|_{i+1}=|x_{i+1}|+\cdots+|x_{n-1}|$. \\
 
 \begin{prop}\label{induced-n-Hom-Lie-rep}
 Let $(V,\rho,\alpha_V)$ be a representation of a multiplicative Hom-Lie Superalgebra $(\mathcal{A},[\cdot,\cdot],\alpha)$ and $\phi\in\wedge^{n-2}\mathcal{A}^*$ satisfying the conditions
 \eqref{NHomLieProduct} and \eqref{NHomLieProduct1}. Then $(V,\rho_\phi,\alpha_V)$ is a representation of the induce $n$-Hom-Lie superalgebra $(\mathcal{A},[\cdot,\cdots,\cdot]_\phi,\alpha)$, where $\rho_\phi$ is defined by \eqref{n-ary-hom-poisson-induce-rep}.
 \end{prop}
 \begin{lem}\label{phi-relation1}
Let $\phi\in\wedge^{n-2}\mathcal{A}^*$ satisfying the conditions
 \eqref{NHomLieProduct} and \eqref{NHomLieProduct1}. Then 
 \begin{equation}\label{condition-phi-crochet-phi}
     \phi([x_1,\cdots,x_n]_\phi,y_1,\cdots,y_{n-3})=0,
 \end{equation}
 for all $x_1,\cdots,x_n,y_1\cdots,y_{n-3}\in\mathcal{H}(\mathcal{A})$, where $[\cdot,\cdots,\cdot]_\phi$ is defined by equation \eqref{nProduct}.
\end{lem}

\begin{proof}
 For all $x_1,\cdots,x_n,y_1\cdots,y_{n-3}\in\mathcal{H}(\mathcal{A})$, we have
 \begin{align*}
     \phi([x_1,\cdots,x_n]_\phi,y_1,\cdots,y_{n-3})&=(-1)^{|X||Y|}\phi(y_1,\cdots,y_{n-3},[x_1,\cdots,x_n]_\phi)\\&=(-1)^{|X||Y|}\displaystyle\sum_{i<j}(-1)^{i+j+1}(-1)^{\gamma^X_{ij}}\phi(x_1,\cdots,\hat{x_i},\cdots,\hat{x_j},\cdots,x_n)\phi(y_1,\cdots,y_{n-3},[x_i,x_j)\\&=(-1)^{|X||Y|}\phi\wedge\delta\phi_Y(X)\\&=0,
 \end{align*}
 since condition \eqref{NHomLieProduct}, where $|X|=|x_1|+\cdots+|x_n|$ and $|Y|=|y_1|+\cdots+|y_{n-3}|$.
\end{proof}

\begin{proof}(Proof of Proposition \ref{induced-n-Hom-Lie-rep})
Let $X=(x_1,\cdots,x_{n-1})\in\mathcal{H}(\mathcal{A})^{n-1}$ and $Y=(y_1,\cdots,y_n)\in\mathcal{H}(\mathcal{A})^{n}$.\\
In the one hand, we have:
\begin{align*}
A&=\rho(\alpha(x_1),\cdots,\alpha(x_{n-1}))\rho(y_1,\cdots,y_{n-1})-(-1)^{|X||Y|^{n-1}}
    \rho(\alpha(y_1),\cdots,\alpha(y_{n-1}))\rho(x_1,\cdots,x_{n-1})\\&=\displaystyle\sum_{i=1}^{n-1}\displaystyle\sum_{j=1}^{n-1}(-1)^{i+j}(-1)^{|x_i||X|_{i+1}+|y_j||Y|_{j+1}^{n-1}}\phi(\alpha(x_1),\cdots,\hat{\alpha(x_i)},\cdots,\alpha(x_{n-1}))\phi(y_1,\cdots,\hat{y_j},\cdots,y_{n-1})\rho(\alpha(x_i))\rho(y_j)\\&-(-1)^{|X||Y|^{n-1}}\displaystyle\sum_{i=1}^{n-1}\displaystyle\sum_{j=1}^{n-1}(-1)^{i+j}(-1)^{|x_i||X|_{i+1}+|y_j||Y|_{j+1}^{n-1}}\phi(\alpha(y_1),\cdots,\hat{\alpha(y_j)},\cdots,\alpha(y_{n-1}))\phi(x_1,\cdots,\hat{x_i},\cdots,x_{n-1})\rho(\alpha(y_j))\rho(x_i)\\&=\displaystyle\sum_{i,j=1}^{n-1}(-1)^{i+j}(-1)^{|x_i||X|_{i+1}+|y_j||Y|_{j+1}^{n-1}}\phi(x_1,\cdots,\hat{x_i},\cdots,x_{n-1})\phi(y_1,\cdots,\hat{y_j},\cdots,y_{n-1})\Big(\rho(\alpha(x_i))\rho(y_j)-(-1)^{|X||Y|^{n-1}}\rho(\alpha(y_j))\rho(x_i)\Big).
 \end{align*}
 On the other hand and using Lemma \ref{phi-relation1}, we have
 \small{\begin{align*}
 B&=\displaystyle\sum_{j=1}^{n-1}(-1)^{|X||Y|^{j-1}}\rho_\phi(\alpha(y_1),\cdots,\alpha(y_{j-1}),[x_1,\cdots,x_{n-1},y_j]_\phi,\alpha(y_{j+1}),\cdots,\alpha(y_{n-1}))\alpha_V\\&=\displaystyle\sum_{j=1}^{n-1}\displaystyle\sum_{i=1}^{j-1}(-1)^{n-i-1}(-1)^{|X||Y|^{j-1}}(-1)^{|y_i|(|Y|_{i+1}+|X|)}\phi(\alpha(y_1),\cdots,\hat{\alpha(y_i)},\cdots,\alpha(y_{j-1}),[x_1,\cdots,x_{n-1},y_j]_\phi,\alpha(y_{j+1}),\cdots,\alpha(y_{n-1}))\rho(\alpha(y_i))\alpha_V\\&+\displaystyle\sum_{j=1}^{n-1}\displaystyle\sum_{i=j+1}^{n-1}(-1)^{n-i-1}(-1)^{|X||Y|^{j-1}}(-1)^{|y_i||Y|_{i+1}}\phi(\alpha(y_1),\cdots,\alpha(y_{j-1}),[x_1,\cdots,x_{n-1},y_j]_\phi,\alpha(y_{j+1}),\cdots,\hat{\alpha(y_i)},\cdots,\alpha(y_{n-1}))\rho(\alpha(y_i))\alpha_V\\&+\displaystyle\sum_{j=1}^{n-1}(-1)^{n-j-1}(-1)^{|X||Y|^{j-1}}(-1)^{|Y|_{j+1}^{n-1}(|y_j|+|X|)}\phi(\alpha(y_1),\cdots,\alpha(y_{j-1}),\alpha(y_{j+1}),\cdots,\alpha(y_{n-1}))\rho([x_1,\cdots,x_{n-1},y_j]_\phi)\alpha_V\\&= \displaystyle\sum_{j=1}^{n-1}(-1)^{n-j-1}(-1)^{|X||Y|^{j-1}}(-1)^{|Y|_{j+1}^{n-1}(|y_j|+|X|)}\phi(\alpha(y_1),\cdots,\alpha(y_{j-1}),\alpha(y_{j+1}),\cdots,\alpha(y_{n-1}))\rho([x_1,\cdots,x_{n-1},y_j]_\phi)\alpha_V\\&=\displaystyle\sum_{i,j=1}^{n-1}(-1)^{i+j}(-1)^{|X||Y|^{j-1}+|x_i||X|_{i+1}}(-1)^{|Y|_{j+1}^{n-1}(|y_j|+|X|)}\phi(\alpha(y_1),\cdots,\alpha(y_{j-1}),\alpha(y_{j+1}),\cdots,\alpha(y_{n-1}))\phi(x_1,\cdots,\hat{x_i},\cdots,x_{n-1})\rho([x_i,y_j])\alpha_V \\+&\displaystyle\sum_{j=1}^{n-1}(-1)^{|X||Y|^{j-1}}(-1)^{|Y|_{j+1}^{n-1}(|y_j|+|X|)}\phi(\alpha(y_1),\cdots,\alpha(y_{j-1}),\alpha(y_{j+1}),\cdots,\alpha(y_{n-1}))\\&\rho\Big(\underbrace{\displaystyle\sum_{1=i<k}^{n-1}(-1)^{i+k+1}(-1)^{(|x_i|+|x_k|)(|X|_{k+1}+|y_j|)+|X|_{i+1}^{k-1}}\phi(x_1,\cdots,\hat{x_i},\cdots,\hat{x_k},\cdots,x_{n-1},y_j)[x_i,x_k]}_{=0\;(By\;eq\; \eqref{NHomLieProduct})}\Big)\alpha_V\\&=\displaystyle\sum_{i,j=1}^{n-1}(-1)^{i+j}(-1)^{|X||Y|^{j-1}+|x_i||X|_{i+1}}(-1)^{|Y|_{j+1}^{n-1}(|y_j|+|X|)}\phi(\alpha(y_1),\cdots,\alpha(y_{j-1}),\alpha(y_{j+1}),\cdots,\alpha(y_{n-1}))\phi(x_1,\cdots,\hat{x_i},\cdots,x_{n-1})\rho([x_i,y_j])\alpha_V
 \end{align*}}
Since the map $\phi$ is even, then we find
 \begin{align*}
 A-B&=\displaystyle\sum_{i,j=1}^{n-1}(-1)^{i+j}(-1)^{|x_i||X|_{i+1}+|y_j||Y|_{j+1}^{n-1}}\phi(x_1,\cdots,\hat{x_i},\cdots,x_{n-1})\phi(y_1,\cdots,\hat{y_j},\cdots,y_{n-1})\Big(\rho(\alpha(x_i))\rho(y_j)-(-1)^{|x_i||y_j|}\rho(\alpha(y_j))\rho(x_i)-\rho[x_i,y_j]\alpha_V\Big)\\&=0.    
 \end{align*}
 The last equality is deduced from equation \eqref{repres-Hom-Lie2} by the fact that $\rho$ is a representation of $(\mathcal{A},[\cdot,\cdot],\alpha)$ which gives that $\rho_\phi$ holds the condition \eqref{repr-n-Hom-Lie1} on $(\mathcal{A},[\cdot,\cdots,\cdot]_\phi,\alpha)$.\\
 
By the same techniques we show that $\rho_\phi$ satisfies the condition \eqref{repr-n-Hom-Lie2}, therefore $(V,\rho_\phi,\alpha_V)$ is a representation of the $n$-Hom-Lie superalgebra $(\mathcal{A},[\cdot,\cdots,\cdot]_\phi,\alpha)$.\\
\end{proof}
\begin{df}
A representation of an $n$-ary Hom-Nambu Poisson superalgebra $(\mathcal{A},[\cdot,\cdots,\cdot],\mu,\alpha)$
 is a quadruple $(V,\rho,\eta,\alpha_V)$ consisting of a $\mathbb{Z}_2$-graded vector space $V$, two aven lineear maps
 $\rho:\Lambda^{n-1}\mathcal{A}\to gl(V)$ and $\eta:\mathcal{A}\to gl(V)$, and an even linear map $\alpha_V:V\to V$ such that
\begin{enumerate}
    \item $(V,\rho,\alpha_V)$ is a representation of the $n$-Hom-Lie Superalgebra $(\mathcal{A},[\cdot,\cdots,\cdot],\alpha)$.
    \item $(V,\eta,\alpha_V)$ is a representation of the Hom-associative superalgebra $(\mathcal{A},\mu,\alpha)$.
    \item The following conditions holds for all $x_1,\cdots,x_{n-1},y,z\in\mathcal{H}(\mathcal{A})$.
   \small{ \begin{eqnarray}
    \eta([x_1,\cdots,x_{n-1},y])\alpha_V&=&\rho(\alpha(x_1),\cdots,\alpha(x_{n-1}))\eta(y)-
    (-1)^{|y||X|}\eta(\alpha(y))\rho(x_1,\cdots,x_{n-1}),\label{cond-rep-hom-poisson1}\\\rho(\alpha(x_1),\cdots,\alpha(x_{n-2}),\mu(y,z))\alpha_V&=&
    (-1)^{|y|(|x_1|+\cdots+|x_{n-2}|)}\eta(\alpha(y))\rho(x_1,\cdots,x_{n-2},z)\label{cond-rep-hom-poisson2}\\
     \nonumber&+&(-1)^{|z|(|x_1|+\cdots+|x_{n-2}|+|y|)}\eta(\alpha(z))\rho(x_1,\cdots,x_{n-2},y),\label{rep-hom-poisson2}
    \end{eqnarray}}
    where $|X|=|x_1|+\cdots+|x_{n-1}|$.
\end{enumerate}
\end{df}
\begin{thm}
Let $(\mathcal{A},[\cdot,\cdots,\cdot],\mu,\alpha)$ be an $n$-ary Hom-Nambu Poisson superalgebra, $\rho:\Lambda^{n-1}\mathcal{A}\to gl(V)$ and $\eta:\mathcal{A}\to gl(V)$ are two even linear maps and $\alpha_V:V\to V$ be an even linear map. Then $(\mathcal{A}\oplus V,[\cdot,\cdots,\cdot]_{\mathcal{A}\oplus V},\mu_{\mathcal{A}\oplus V},\alpha+\alpha_V)$ is an $n$-ary Hom-Nambu Poisson superalgebra if and only if $(V,\rho,\eta,\alpha_V)$ is a representation of $(\mathcal{A},[\cdot,\cdots,\cdot],\mu,\alpha)$, where $\mu_{\mathcal{A}\oplus V}$,  $(\alpha+\alpha_V)$ and $[\cdot,\cdots,\cdot]_{\mathcal{A}\oplus V}$ are defined respectively by the eqs. \eqref{Hom-ass-direct-sum1}, \eqref{Hom-ass-direct-sum2} and \eqref{crochet-direct-sum-n-hom-lie-superalgebras}.
\end{thm}
\begin{proof}
 By Proposition \ref{direct-sum-n-hom-lie-superalgebras}, $(\mathcal{A}\oplus V,[\cdot,\cdots,\cdot]_{\mathcal{A}\oplus V},\alpha+\alpha_V)$ is an $n$-Hom-Lie superalgebras if and only if $(V,\rho,\alpha_V)$ is a representation of $(\mathcal{A},[\cdot,\cdots,\cdot],\alpha)$ and Proposition \ref{sumdirecthomassoc} gives that $(\mathcal{A}\oplus V,\mu_{\mathcal{A}\oplus V},\alpha+\alpha_V)$ is a commutative Hom-associative superalgebra if and only if $(V,\eta,\alpha_V)$ is a representation of $(\mathcal{A},\mu,\alpha)$.\\
 Moreover, it is easy to show that the identity \eqref{n-aryHomLeibniz} is satisfied on $(\mathcal{A}\oplus V,[\cdot,\cdots,\cdot]_{\mathcal{A}\oplus V},\mu_{\mathcal{A}\oplus V},\alpha+\alpha_V)$ if and only if conditions \eqref{cond-rep-hom-poisson1}  and \eqref{cond-rep-hom-poisson2} are satisfied which gives the proof of the theorem.
 \end{proof}

\begin{thm}\label{induced-repr-n-Poisson}
Let $(V,\rho,\eta,\alpha_V)$ be a representation of a multiplicative Hom-Poisson superalgebra $(\mathcal{A},[\cdot,\cdot],\mu,\alpha)$ and $\phi\in\wedge^{n-2}\mathcal{A}^*$ satisfying the conditions
 \eqref{NHomLieProduct} and \eqref{NHomLieProduct1}. Then $(V,\rho_\phi,\eta,\alpha_V)$ is a representation of the induce Hom-Poisson superalgebra $(\mathcal{A},[\cdot,\cdots,\cdot]_\phi,\mu,\alpha)$, where $\rho_\phi$ is defined by \eqref{n-ary-hom-poisson-induce-rep}.\\
\end{thm}

\begin{proof} 
By Proposition \ref{induced-n-Hom-Lie-rep}, the triple $(V,\rho_\phi,\alpha_V)$ is a representation of the induced $n$-Hom Lie superalgebra $(\mathcal{A},[\cdot,\cdots,\cdot]_\phi,\alpha)$.\\

It is obvious that $(V,\eta,\alpha)$ is a representation of the Hom-associative superalgebra $(\mathcal{A},\mu,\alpha)$. It remains to show that $\rho_\phi$ satisfies the conditions \eqref{cond-rep-hom-poisson1} and  \eqref{cond-rep-hom-poisson2}.\\

Let $x_1,\cdots,x_{n-1},y\in\mathcal{H}(\mathcal{A})$, we have

\begin{align*}
\eta([x_1,\cdots,x_{n-1},y]_\phi)\alpha_V&=\eta\Big(\underbrace{\displaystyle\sum_{1=i<j}^{n-1}(-1)^{i+j+1}(-1)^{\gamma_{ij}^X}(-1)^{|y|(|x_i|+|x_j|)}\phi(x_1,\cdots,\hat{x_i},\cdots,\hat{x_j},\cdots,x_{n-1},y)[x_i,x_j]}_{=0\;(By\;eq\; \eqref{NHomLieProduct})}\Big)\alpha_V \\&+\displaystyle\sum_{i=1}^{n-1}(-1)^{i+n+1}(-1)^{|x_i|(|X|_{i+1}+|y|)}\phi(x_1,\cdots,\hat{x_i},\cdots,x_{n-1})\eta([x_i,y])\alpha_V, 
\end{align*}

and 

\begin{align*}
 \rho_\phi(\alpha(x_1),\cdots,\alpha(x_{n-1}))\eta(y)&-(-1)^{|y||X|}\eta(\alpha(y))\rho_\phi(x_1,\cdots,x_{n-1})\\&=\displaystyle\sum_{i=1}^{n-1}(-1)^{n-i-1}(-1)^{|x_i||X|_{i+1}}\phi(\alpha(x_1),\cdots,\hat{\alpha(x_i)},\cdots,\alpha(x_{n-1}))\rho(x_i)\eta(y)\\&-(-1)^{|y||X|}\eta(\alpha(y))\displaystyle\sum_{i=1}^{n-1}(-1)^{n-i-1}(-1)^{|x_i||X|_{i+1}}\phi(x_1,\cdots,\hat{x_i},\cdots,x_{n-1})\rho(x_i)\\&=\displaystyle\sum_{i=1}^{n-1}(-1)^{n+i+1}(-1)^{|x_i||X|_{i+1}}\phi(x_1,\cdots,\hat{x_i},\cdots,x_{n-1})\Big(\rho(\alpha(x_i))\eta(y)-(-1)^{|y||x_i|}\eta(\alpha(y))\rho(x_i)\Big).  
\end{align*}
By eq. \eqref{repr-Hom-Poisson1},  $\eta([x_i,y])\alpha_V=\rho(\alpha(x_i))\eta(y)-(-1)^{|y||x_i|}\eta(\alpha(y))\rho(x_i)$ which gives that $\eta([x_1,\cdots,x_{n-1},y]_\phi)\alpha_V=\rho_\phi(\alpha(x_1),\cdots,\alpha(x_{n-1}))\eta(y)-(-1)^{|y||X|}\eta(\alpha(y))\rho_\phi(x_1,\cdots,x_{n-1})$ and the condition \eqref{cond-rep-hom-poisson1} is verified.\\

By the same way, we show that condition \eqref{cond-rep-hom-poisson2} is verified. The theorem is proved. 
\end{proof}

\end{document}